\documentclass[a4paper,reqno]{amsart}
\usepackage{amsmath, amsthm, amsfonts, amssymb, mathrsfs}
\usepackage{enumitem}
\usepackage{graphicx}
\usepackage[update,prepend]{epstopdf} 
\usepackage{color}
\usepackage{subfigure}

\newcommand{\beq}{\begin{equation}}
\newcommand{\eeq}{\end{equation}}
\newcommand{\be}{\begin{equation*}}
\newcommand{\ee}{\end{equation*}}
\newcommand{\bmat}{\begin{pmatrix}}
\newcommand{\emat}{\end{pmatrix}}

\newcommand{\N}{\mathbb{N}}

\newcommand{\R}{\mathbb{R}}
\newcommand{\C}{\mathbb{C}}

\newcommand\I{\mathrm{i}}
\newcommand\re{\mathrm{Re}}
\newcommand\im{\mathrm{Im}}

\newcommand{\rd}{\mathrm{d}}

\newcommand\dom{\mathcal D}
\newcommand\ran{\mathcal R}
\newcommand\cP{\mathcal P}

\newcommand\cS{\mathcal S}
\newcommand\cT{\mathcal T}

\newcommand\cA{\mathcal A}
\newcommand\cB{\mathcal B}

\newcommand\cE{\mathcal E}

\newcommand\lbar\overline

\newcommand\eps\varepsilon
\renewcommand\epsilon\varepsilon
\renewcommand\phi\varphi
\renewcommand\rho\varrho
\newcommand\al\alpha
\newcommand\lm\lambda

\newcommand\ds\displaystyle

\newcommand\p\partial

\newcommand{\diag}{{\rm diag} }

\newcommand{\To}{\Longrightarrow}

\newcommand{\tolong}{\longrightarrow}

\newcommand{\gnr}{\stackrel{gnr}{\rightarrow}}
\newcommand{\gsr}{\stackrel{gsr}{\rightarrow}}
\newcommand{\s}{\stackrel{s}{\rightarrow}}

\newcommand{\slong}{\stackrel{s}{\longrightarrow}}

\newcommand{\dist}{{\rm dist}}

\begin{document}

\theoremstyle{plain}
\newtheorem{theorem}{Theorem}[section]
\newtheorem{lemma}[theorem]{Lemma}
\newtheorem{criterion}[theorem]{Criterion}
\newtheorem{prop}[theorem]{Proposition}
\newtheorem{corollary}[theorem]{Corollary}
\renewcommand{\proofname}{Proof}
\theoremstyle{definition} 
\newtheorem{definition}[theorem]{Definition}
\newtheorem{example}[theorem]{Example}
\newtheorem{rem}[theorem]{Remark}
\newtheorem{notation}[theorem]{Notation}
\newtheorem{ass}[theorem]{Assumption}

\title{Convergence of sequences of linear operators and their spectra}

\author{Sabine B\"ogli}
\address[S.\ B\"ogli]{
Department of Mathematics,
Imperial College London,
Huxley Building,
180 Queen's Gate, London SW7 2AZ, UK}
\email{s.boegli@imperial.ac.uk}

\subjclass[2010]{47A10, 47A55, 47A58, 47B07}

\keywords{eigenvalue approximation, spectral exactness, generalized resolvent convergence, discrete compactness}

\date{
\today
}

\begin{abstract}
We establish spectral convergence results of approximations of unbounded non-selfadjoint linear operators with compact resolvents by operators that converge in generalized strong resolvent sense.
The aim is to establish general assumptions that ensure spectral exactness, i.e.\ that every true eigenvalue is approximated and no spurious eigenvalues occur.
A main ingredient is the discrete compactness of the sequence of resolvents of the approximating operators.
We establish sufficient conditions and perturbation results for strong convergence and for discrete compactness of the resolvents.
\vspace{-0.7mm}
\end{abstract}

\maketitle

\section{Introduction}
The spectra of linear operators $T$, e.g.\ describing the time evolution of a physical system, are usually not known analytically 
and need to be computed numerically by approximating the operators and determining the eigenvalues of simpler operators. 
However, it is well-known that spectral computations may lead to \emph{spectral pollution}, i.e.\ to 
numerical artefacts which do not belong to the spectrum of $T$, so-called \emph{spurious eigenvalues}.
Vice versa, not every eigenvalue or spectral point of $T$ may be approximated; an approximation $(T_n)_{n\in\N}$ is called \emph{spectrally inclusive}
if this phenomenon does not occur. 
If spectral inclusion prevails and no spectral pollution occurs, then $(T_n)_{n\in\N}$ is said to be a \emph{spectrally exact} approximation of $T$.

{The existing spectral exactness results in the literature are restricted either  to 
bounded operators or 
to particular classes of differential operators, or
they are only local spectral exactness results, e.g.\  for spectral gaps of selfadjoint operators.
On the other hand, many important applications in physics such as linear stability problems in fluid mechanics,
magnetohydrodynamics, or elasticity theory require reliable knowledge on the spectra of unbounded non-selfadjoint linear operators.

The present paper aims at filling this gap.
The novelty of the results established here lies in
1) their far-reaching generality covering wide classes of unbounded non-self\-adjoint linear operators;
2) their simultaneous applicability to different approximation schemes such as the Galerkin (finite section) method 
and the domain truncation method;
3) their global nature which yields spectral exactness in the \emph{entire} complex plane; and
4) a comprehensive analysis of necessary conditions and perturbation results for spectral exactness.
We present applications to interval truncation of singular $2\times 2$ differential operator matrices, to domain truncation of magnetic Schr\"odinger operators with complex-valued potentials on $\R^d$, and to the Galerkin method for operators of block-diagonally dominant form.
The efficacy of the results of this paper is substantiated by the fact that they form the basis to establish spectral exactness
for all presently considered magnetohydrodynamic (MHD) dynamo models in astrophysics
for which \emph{no} analytic results whatsoever have been available up to now (see~\cite{dynamopaper},  based on~\cite[Chapter~3]{boegli-phd}).}

The first main theorem (Theorem~\ref{mainthmforgsr}) is the following global spectral convergence result:
If $(T_n)_{n\in\N}$ converges in generalized strong resolvent (gsr) sense to~$T$ and the resolvents are compact and form a discretely compact sequence,
then the approximation is spectrally exact.
In the second main result  (Theorem~\ref{thmgsrimpliesgnrbasic}) we prove that under the additional assumptions that the operators act in Hilbert spaces and $(T_n^*)_{n\in\N}$ converges to $T^*$ in gsr-sense, then the resolvents converge even in operator norm. 
The third group of important results comprise additive perturbation results:
For a sequence $(S_n)_{n\in\N}$ of relatively bounded perturbations $S_n$ of $T_n$ with (uniform) relative bound~$<1$,
we establish perturbation results for gsr-convergence (Theorem~\ref{perturbgsr})
and discrete resolvent compactness (Theorem~\ref{perturbdiscrcompres}).
A fourth group of results guarantee gsr-convergence and discrete resolvent compactness of a sequence of block operator matrices
by means of easily verifiable assumptions that are formulated in terms of the matrix entries. 
First we prove results for unbounded \emph{finite} operator matrices 
and then for \emph{infinite} matrices (Theorems~\ref{gsrdiag} and~\ref{discrcompresdiag}).

The notions of spectral inclusion and spectral exactness w{ere} introduced by Bailey et al.\ \cite{Bailey-1993} for regular approximations of singular selfadjoint Sturm-Liouville problems via interval truncation.
{They were further studied, in particular, by Brown and Marletta for the  domain truncation procedure of non-selfadjoint differential operators~\cite{Brown-Marletta-2001,Brown-Marletta-2003,Brown-Marletta-2004}. 
%
%
The notion of  \emph{generalized norm/strong resolvent convergence}
developed in this paper
for approximations of unbounded linear operators
 is closely related to norm or strong resolvent convergence studied by Kato~\cite[Sections~IV.2,VIII.1]{kato},
Reed-Simon~\cite[Theorems~VIII.23-25]{reedsimon1} and Weidmann~\cite[Section~9.3]{weid1};
in the latter two, only selfadjoint operators were considered.
In general, the approximating operators $T_n$ 
cannot be chosen to act in the same space as $T$,
so we 
compare the projected resolvents $(T_n-\lm)^{-1}P_n$ and $(T-\lm)^{-1}P$ in a common larger space. 
Note that the meaning of ``generalized'' used in this paper 
is different from Kato's \emph{generalized convergence}
(meaning resolvent convergence)
where it indicates that the operators are unbounded.
The spectral exactness result  (Theorem~\ref{mainthmforgsr}) relies on gsr-convergence, however we are also interested 
in generalized \emph{norm} resolvent convergence (see Theorem~\ref{thmgsrimpliesgnrbasic}) since the latter is used to prove convergence of pseudospectra in Hausdorff metric (see~\cite[Theorem~2.1]{boegli-siegl-2014}).


To conclude spectral exactness, it is not enough to assume that the operators $T$ and $T_n$, $n\in\N$, have compact resolvents and converge in gsr-sense.
In fact, even in the selfadjoint case, if the operator $T$ is unbounded below and above, then the Galerkin method may produce spurious eigenvalues \emph{anywhere} on the real line (see~\cite[Theorem~2.1]{levitinshargorodsky}). 
We prove spectral exactness under the additional assumption that the sequence $((T_n-\lm)^{-1})_{n\in\N}$ is \emph{discretely compact}.
The latter notion was introduced by Stummel who  established a spectral convergence theory for bounded operators in \cite{stummel1,stummel2}.
Similar result were obtained by  Anselone-Palmer and Osborn~\cite{anselonepalmer,osborn} for the closely related notion of \emph{collectively compact} sets of bounded operators.
%
For relations between the various results and notions,
we refer to Chatelin's monograph~\cite{chatelin} (see, in particular, Sections~3.1-3.6,~5.1-5.5). 

We mention that if the assumptions of Theorem~\ref{mainthmforgsr} are not satisfied (in particular if essential spectrum is present), a different approach is to establish \emph{local} spectral exactness results, 
i.e.\ to identify regions in the complex plane where no spectral pollution occurs or to find enclosures for true eigenvalues.
This was done for the Galerkin approximation of selfadjoint operators by means of 
{\emph{higher order relative spectra} (introduced by Davies in~\cite{davies_higher_order}, see also the  comprehensive overviews by Shargorodsky et al.\ \cite{shargorodsky-higher-order-relative-spectra,levitinshargorodsky})}
and the closely related methods of Davies-Plum~\cite{davies-plum-2004} and Mertins-Zimmermann~\cite{zimmermann-mertins}.
For non-selfadjoint operators, we prove local spectral exactness results in terms of the \emph{region of boundedness} (see Theorem~\ref{thmlocal}).
An alternative but computationally very expensive method 
to obtain reliable information on isolated eigenvalues uses \emph{interval arithmetic} which yields eigenvalue enclosures with absolute certainty (see e.g.\ Brown et al.\ \cite{blmtw2} and the references therein).

This paper is organized as follows.
In Section~\ref{sectionspectralconv} we prove convergence results for operators and their spectra
under the assumptions that the operators have compact and discretely compact resolvents and converge in gsr-sense.
In Section~\ref{sectionsuffcondforstrongresconv} and Section~\ref{sectionsuffcondforcompres} we derive sufficient conditions and perturbation results for gsr-convergence and for discrete resolvent compactness, respectively.
In both sections 
we complement the general theorems by results for finite and for infinite unbounded operator matrices.
Applications to the domain truncation method for singular differential operators (and operator matrices) and to the Galerkin method are given in Section~\ref{sectionapplications}.


We use the following notation.
The norm of a normed space $E$ is denoted by~$\|\cdot \|_E$.
The convergence in $E$, i.e.\ $\|x_n-x\|_E\to 0$, is written as $x_n\to x$.
In a Hilbert space~$H$ the scalar product is $\langle\cdot,\cdot\rangle_H$.
Weak convergence in $H$, i.e.\ $\langle x_n,z\rangle_H\to \langle x,z\rangle _H$ for all $z\in H$, is denoted by $x_n\stackrel{w}{\to} x$.
For two normed vector spaces $D$ and $E$ we denote by $L(D,E)$ the space of all bounded linear operators from $D$ to $E$; we write $L(E)$ if $D=E$. 
Analogously, the space of all closed operators in $E$ is denoted by $C(E)$.
The spectrum, point spectrum, approximate point spectrum and resolvent set 
of a linear operator $T$ are denoted by $\sigma(T)$, $\sigma_p(T)$, $\sigma_{\rm app}(T)$ and $\rho(T)$, 
respectively, and the Hilbert space adjoint operator of $T$ is~$T^*$.
For an operator $T\in C(E)$ the graph norm is $\|\cdot\|_T:=\|\cdot\|_E+\|T\cdot\|_E$; then $(\dom(T),\|\cdot\|_T)$ is a Banach space.
For bounded linear operators we write $T_n\to T$ and $T_n\s T$ for norm and strong convergence in $L(D,E)$.
An identity operator in a Banach or Hilbert space is denoted by $I$; scalar multiples $\lm I$ are written as $\lm$. Analogously, the 
operator of multiplication with a function $m$ in some $L^2$-space is also denoted by $m$.
Given $\lm\in\C$ and a subset $\Omega\subset\C$, the distance of $\lm$ to $\Omega$ is ${\rm dist}(\lm,\Omega):=\inf_{z\in\C}|z-\lm|$, 
and $B_r(\lm):=\{z\in\C:\,|z-\lm|<r\}$ for $r>0$.
Let $\N:=\{1,2,3,\dots\}$, in particular, $0\notin\N$. 
Finally, for a subset $I\subset\N$, we denote by $\# I$ the number of elements in $I$.

\section{Convergence of operators and their spectra}\label{sectionspectralconv}

In this section we establish convergence results for operators acting in different spaces and spectral convergence results.
In Subsection~\ref{subsectionmainresults} are the main convergence results. 
Before we prove these results  in Subsection~\ref{subsectionproofs}, we recall
the notions of discretely compact sequences or collectively compact sets of bounded operators 
and we analyze their effect on strong or norm operator convergence (see Subsection~\ref{subsectionbounded}).

In the following, we assume  that $E_0$ is a Banach space and $E, \,E_n\subset E_0, \,n\in\N,$ are closed subspaces.
Let $P:E_0\to E$ and $P_n:E_0\to E_n$, $n\in\N$, be projections onto the respective subspaces converging strongly, $P_n\s P$.
Throughout, in results for Hilbert spaces $H_0:=E_0,\,H:=E,\,H_n:=E_n$, $n\in\N$, we assume that $P$, $P_n$, $n\in\N$, are the \emph{orthogonal} projections onto the respective subspaces.

\subsection{Main convergence results for unbounded linear operators and their spectra}\label{subsectionmainresults}
%


The following definition of \emph{generalized} strong and norm resolvent convergence is due to Weidmann~\cite[Section~9.3]{weid1},
and the region of boundedness was introduced by Kato~\cite[Section~VIII.1]{kato}.

\begin{definition}\label{defresconv}
Let $T\in C(E)$ and $T_n\in C(E_n)$, $n\in\N$. 
\begin{enumerate}[label=\rm{\roman{*})}]
\item The sequence $(T_n)_{n\in\N}$ is said to \emph{converge in generalized strong resolvent sense}\index{convergence of operator sequence!generalized strong resolvent convergence} to $T$, 
$T_n\gsr T$,\index{$T_n\gsr T$} if there exist $n_0\in\N$ and $\lm\in\underset{n\geq n_0}{\bigcap}\rho(T_n)\cap\rho(T)$ with\vspace{-2mm}
$$(T_n-\lm)^{-1}P_n\slong (T-\lm)^{-1}P, \quad n\to\infty.$$
\item The sequence $(T_n)_{n\in\N}$ is said to \emph{converge in generalized norm resolvent sense}\index{convergence of operator sequence!generalized norm resolvent convergence} to $T$, 
$T_n\gnr T$,\index{$T_n\gnr T$} if there exist $n_0\in\N$ and $\lm\in\underset{n\geq n_0}{\bigcap}\rho(T_n)\cap\rho(T)$ with\vspace{-2mm}
$$(T_n-\lm)^{-1}P_n\tolong (T-\lm)^{-1}P, \quad n\to\infty.$$
\item The \emph{region of boundedness}\index{region of boundedness!of $(T_n)_{n\in\N}$} of the sequence $(T_n)_{n\in\N}$ is defined as \index{$\Delta_b\left((T_n)_{n\in\N}\right)$}
\begin{align*}&\Delta_b\left((T_n)_{n\in\N}\right):=\left\{\lm\in\C:\,\exists\,n_0\in\N\,\text{with}\begin{array}{l}\lm\in\rho(T_n),\,n\geq n_0,\\ \left(\|(T_n-\lm)^{-1}\|\right)_{n\geq n_0} \text{bounded}\end{array}\right\}.\end{align*}
\end{enumerate}
\end{definition}

The following notions of spectral inclusion and spectral exactness were introduced by Bailey et al.\ in~\cite{Bailey-1993}. We will also use local versions of these notions.
\begin{definition}\label{defspectralexactness}
Let $T\in C(E)$ and $T_n\in C(E_n)$, $n\in\N$. 
\begin{enumerate}[label=\rm{\roman{*})}]
\item The approximation $(T_n)_{n\in\N}$ of $T$ is called \emph{spectrally inclusive}\index{spectral approximation!spectral inclusion} if for every $\lm\in\sigma(T)$ there exists a sequence $(\lm_n)_{n\in\N}$ of elements $\lm_n\in\sigma(T_n)$, $n\in\N$, with $\lm_n\to\lm$.
\item An element $\lm\in\C$ for which there exist an infinite subset $I\subset\N$ and $\lm_n\in\sigma(T_n)$, $n\in I$, with $\lm_n\to\lm$ but $\lm\notin\sigma(T)$ is called 
\emph{spurious eigenvalue}\index{eigenvalue!spurious}. The occurrence of such a point is known as \emph{spectral pollution}\index{spectral approximation!spectral pollution}.
\item  The approximation $(T_n)_{n\in\N}$ of $T$ is called \emph{spectrally exact}\index{spectral approximation!spectral exactness} if it is spectrally inclusive and no spectral pollution occurs. 
\end{enumerate}
\end{definition}

The following result yields local spectral exactness  in the region of boundedness.
\begin{theorem}\label{thmlocal}
Let $T\in C(E)$ and $T_n\in C(E_n)$, $n\in\N$.
Suppose that $T_n\gsr T$.
\begin{enumerate}
\item[\rm i)] For each $\lm\in\sigma(T)$ such that for some $\eps>0$ we have 
\beq B_{\eps}(\lm)\backslash\{\lm\}\subset \Delta_b\left((T_n)_{n\in\N}\right)\cap\rho(T),\label{eqneighbourhoodinregionofbdd}\eeq  
there exist $\lm_n\in\sigma(T_n), \,n\in\N,$ with $\lm_n\to\lm$ as $n\to\infty$.
\item[\rm ii)] No spectral pollution occurs in $\Delta_b\left((T_n)_{n\in\N}\right)$.
\end{enumerate}
\end{theorem}

Next we generalize results that are known for the special case $E_0=E=E_n$  and resolvent convergence to \emph{generalized} resolvent convergence.
Claim~i) is a generalization of \cite[Theorem~IV.2.25]{kato} and  \cite[Satz~9.26~b)]{weid1}; the latter result only applies to selfadjoint operators.
Claim~ii) is a generalization of~\cite[Satz 9.24 a)]{weid1}.

\begin{theorem}\label{thmgeneralizedresconv}
Let $T\in C(E)$ and $T_n\in C(E_n)$, $n\in\N$.
\begin{enumerate}
\item[\rm i)]
If $T_n\gnr T$, then no spectral pollution occurs.

\item[\rm ii)]
Let $E_0$, $E$, $E_n$, $n\in\N$, be Hilbert spaces, and let $T$, $T_n$, $n\in\N$, be selfadjoint with $T_n\gsr T$.
Then \beq\sigma(T)\subset\C\backslash \Delta_b\big((T_n)_{n\in\N}\big)=\big\{\lm\in\C:\,\liminf_{n\to\infty}\,{\rm dist}(\lm,\sigma(T_n))=0\big\},\label{eq.Deltab.selfadjoint}\eeq
and hence $(T_n)_{n\in\N}$ is a spectrally inclusive approximation of $T$.
\end{enumerate}
\end{theorem}

To formulate the main results (Theorems~\ref{mainthmforgsr} and~\ref{thmgsrimpliesgnrbasic} below), we use Stummel's notions of
discrete compactness of a sequence of bounded operators (see~\cite[Definition~3.1.(k)]{stummel1}).

\begin{definition}\label{defdiscreteconvofop}
Let $D_n$, $n\in\N$, be arbitrary Banach spaces and
$A_n\in L(D_n, E_n)$, $n\in\N$. 
The sequence  $(A_n)_{n\in\N}$ is said to be \emph{discretely compact}\index{compactness!discretely compact operator sequence} if for each infinite subset $I\subset\N$ and each bounded sequence of elements $x_n\in D_n, \,n\in I,$ there 
exist $x\in E$ and an infinite subset $\widetilde I\subset I$ so that $\|x_n-x\|_{E_0}\to 0$ as $n\in \widetilde I$, $n\to\infty$. 
\end{definition}

The following theorem is  the main spectral convergence result of this section.
\begin{theorem}\label{mainthmforgsr}
Let $T\in C(E)$ and $T_n\in C(E_n)$, $n\in\N$.
Assume that there exists an element
$\lm_0\in\underset{n\in\N}{\bigcap}\rho(T_n)\cap\rho(T)$ such that $(T-\lm_0)^{-1}$, $(T_n-\lm_0)^{-1},\,n\in\N,$ are compact operators
and the sequence $\left((T_n-\lm_0)^{-1}\right)_{n\in\N}$ is discretely compact.
If $T_n\gsr T$, then the following hold:
\begin{enumerate}[label=\rm{\roman{*})}]
\item The region of boundedness coincides with the resolvent set of $T$, $$\Delta_b\left((T_n)_{n\in\N}\right)=\rho(T),$$ 
and, for any $\lm\in\rho(T)$, 
\begin{equation}\label{eq.gsr.atalllm}
(T_n-\lm)^{-1}P_n\slong (T-\lm)^{-1}P, \quad n\to\infty.
\end{equation}
\item The sequence $(T_n)_{n\in\N}$ is a spectrally exact approximation of $T$.
More precisely, no spectral pollution occurs, and if $\lm \in \C$ is an eigenvalue of $T$ of algebraic multiplicity $m$,
then, for $n$ large enough,  $T_n$ has exactly $m$ eigenvalues {\rm(}repeated according to their algebraic multiplicities{\rm)} in a neighbourhood of $\lm$ which converge to $\lm$ as $n\to \infty$ and the corresponding normalized elements of the algebraic eigenspaces converge {\rm(}with respect to $\|\cdot\|_{E_0}${\rm)}.
\end{enumerate}
\end{theorem}

Now we assume that the underlying spaces are Hilbert spaces. 
We establish sufficient conditions guaranteeing that generalized strong resolvent convergence implies generalized norm resolvent convergence.

\begin{theorem}\label{thmgsrimpliesgnrbasic}
Let $T\in C(H)$ and $T_n\in C(H_n)$, $n\in\N$.
Assume that there exists an element
$\lm_0\in\underset{n\in\N}{\bigcap}\rho(T_n)\cap\rho(T)$ such that  $(T_n-\lm_0)^{-1},\,n\in\N,$ are compact operators
and the sequence $\left((T_n-\lm_0)^{-1}\right)_{n\in\N}$ is discretely compact.
If $$(T_n-\lm_0)^{-1}P_n\slong(T-\lm_0)^{-1}P, \quad (T_n^*-\lbar{\lm_0})^{-1}P_n\slong (T^*-\lbar{\lm_0})^{-1}P,\quad n\to\infty,$$ 
then, for every $\lm\in\rho(T)$, the operator $(T-\lm)^{-1}$ is compact and $$(T_n-\lm)^{-1}P_n\tolong (T-\lm)^{-1}P, \quad n\to\infty.$$
\end{theorem}

\subsection{Convergence and compactness concepts for bounded operators}\label{subsectionbounded}
In this subsection we study discretely compact operator sequences and the effect of this notion on strong operator convergence.

First we prove multiplicative and additive perturbation results on discrete compactness. 
Denote by $D_n$, $n\in\N$, arbitrary Banach spaces.

\begin{lemma}\label{lemmapertdisccomp}
 \begin{enumerate}[label=\rm{\roman{*})}]
  \item Let $A_n\in L(E_n)$, $B_n\in L(D_n,E_n)$, $n\in\N$. If $(A_n)_{n\in\N}$ is discretely compact and $(B_n)_{n\in\N}$ is a bounded sequence, then $(A_nB_n)_{n\in\N}$ is discretely compact.
  \item Let  $A_n\in L(E_n)$, $B_n\in L(D_n,E_n)$, $n\in\N$. If $(B_n)_{n\in\N}$ is discretely compact and there exists $A\in L(E)$ with $A_nP_n\s AP$, then $(A_nB_n)_{n\in\N}$ is discretely compact.
  \item For $j=1,\dots,k$, let $A_n^{(j)}\in L(D_n,E_n)$, $n\in\N$. If the sequences $\big(A_n^{(j)}\big)_{n\in\N}$, $j=1,\dots,k$, are discretely compact, then so is \vspace{-1mm}$$\bigg(\sum_{j=1}^k A_n^{(j)}\bigg)_{n\in\N}.$$
 \end{enumerate}
\end{lemma}
\begin{proof}
i) Let $I\subset\N$ be an infinite subset, $M>0$ and $x_n\in L(D_n)$, $n\in I$, with $\|x_n\|_{D_n}\leq M$ for all  $n\in I$. Define $C:=\sup_{n\in I}\|B_n\|$.
Then $\|B_nx_n\|_{E_n}\leq CM$, $n\in I$. Now the discrete compactness of $(A_n)_{n\in\N}$ (see Definition~\ref{defdiscreteconvofop}) implies that a subsequence of $(A_nB_nx_n)_{n\in I}$ is convergent in $E_0$ with limit in $E$.
So $(A_nB_n)_{n\in\N}$ is discretely compact.

ii) Let $I\subset\N$ be an infinite subset, $M>0$ and $x_n\in L(D_n)$, $n\in I$, with $\|x_n\|_{D_n}\leq M$,  $n\in I$.
The discrete compactness of $(B_n)_{n\in\N}$ implies that there exist $y\in E$ and an infinite subset $\widetilde I\subset I$ such that $\|B_nx_n-y\|_{E_0}\to 0$ as $n\in \widetilde I$, $n\to\infty$.
By the assumptions, we have $A_nP_n\s AP$, so the element $z:=Ay\in E$ satisfies $\|A_nB_nx_n-z\|_{E_0}\to 0$. Hence $(A_nB_n)_{n\in \N}$ is discretely compact.

iii) Let $I\subset\N$ be an infinite subset, $M>0$ and $x_n\in L(D_n)$, $n\in I$, with  $\|x_n\|_{D_n}\leq M$, $n\in I$. 
By the discrete compactness of $\big(A_n^{(1)}\big)_{n\in\N}$, there exists an infinite subset  $I^{(1)}\subset I$ such that the sequence $\big(A_n^{(1)}x_n\big)_{n\in I^{(1)}}$ is  convergent in $E_0$ with limit in~$E$.
Now, inductively for every $j=2,\dots,k$, the discrete compactness of $\big(A_n^{(j)}\big)_{n\in\N}$ implies that there exists an infinite subset $I^{(j)}\subset I^{(j-1)}$ such that $\big(A_n^{(j)}x_n\big)_{n\in I^{(j)}}$ is  convergent in $E_0$ with limit in~$E$.
Therefore \vspace{-1mm}
$$\bigg(\sum_{j=1}^k A_n^{(j)}x_n\bigg)_{n\in I^{(k)}}\vspace{-1mm}$$ is convergent in $E_0$ with limit in~$E$. 
\end{proof}

Analogously to the result that the limit of a sequence of compact operators is compact (see e.g.\ \cite[Theorem III.4.7]{kato}), one can show the following result for discrete compactness.
An application of Proposition~\ref{propapproxdiscrcomp} is given in Theorem~\ref{discrcompresdiag} where infinite diagonal operator matrices are approximated by $k\times k$ matrices. 
\begin{prop}
Let $E_0=E$. 
For each $n\in\N$ let $A_n, \,A_n^{(k)}\in L(D_n,E_n), \,k\in\N,$ with 
\beq \sup_{n\in\N}\big\|A_n^{(k)}- A_n\big\|\tolong 0, \quad k\to\infty.\label{eq.unif.conv}\eeq
If all sequences $\big(A_n^{(k)}\big)_{n\in\N}, \,k\in\N,$ are discretely compact, then so is~$(A_n)_{n\in\N}$.
\label{propapproxdiscrcomp}
\end{prop}

\begin{proof}
Consider an infinite subset $I\subset\N$ and a bounded sequence of elements $x_n\in D_n, \,n\in I,$ i.e.\ there exists $M>0$ such that $\|x_n\|_{D_n}\leq M, \,n\in I$. 
 We show the existence of an infinite subset $\widetilde I\subset I$ such that for all $\eps>0$ there exists $N_{\eps}\in\N$ with 
 $$\|A_nx_n-A_mx_m\|_{E}<\eps, \quad n,m\in \widetilde I, \,n,m\geq N_{\eps}.$$
 Then the claim follows from the completeness of $E$.
 
 The sequence $\big(A_n^{(1)}\big)_{n\in\N}$ is discretely compact by the assumptions, hence there exists an infinite subset $I^{(1)}\subset I$ such that 
 $\big(A_n^{(1)}x_n\big)_{n\in I^{(1)}}$ is convergent in $E$. 
Inductively, for each $k\geq 2$, we find  an infinite subset $I^{(k)}\subset I^{(k-1)}$ such that  
 $\big(A_n^{(k)}x_n\big)_{n\in I^{(k)}}$ is convergent in $E$. Therefore, there exists an increasing sequence $(N^{(k)})_{k\in\N}\subset\N$ such that $N^{(k)}\in I^{(k)}, \,k\in\N,$ and
 $$\big\|A_n^{(k)}x_n-A_m^{(k)}x_m\big\|_{E}<\frac{1}{k}, \quad n,m\in I^{(k)},\quad n,m\geq N^{(k)}.$$
 %
We define $\widetilde I:=\big\{N^{(k)}:\,k\in\N\big\}.$
 Let $\eps>0$ be fixed.  
By the assumption~\eqref{eq.unif.conv}, we find $K_{\eps}\in\N$ such that $K_{\eps}\geq\frac{3}{\eps}$ and  
 $$\big\|A_n^{(k)}-A_n\big\|<\frac{\eps}{3M}, \quad n\in\N, \quad k\geq K_{\eps}.$$
 Altogether, for $k\geq K_{\eps}$ (which yields $\frac{1}{k}\leq \frac{1}{K_{\eps}}\leq\frac {\eps}{3}$) and $l\geq k$, the elements $n=N^{(k)}$, $m=N^{(l)}\in\widetilde I$ satisfy $n,m\geq N^{(K_{\eps})}=:N_{\eps}$ and
 \begin{align*}
 &\left\|A_{n}x_{n}-A_{m}x_{m}\right\|_{E}\\
 &\leq \big\|A_{n}-A_{n}^{(k)}\big\|\,\|x_{n}\|_{E_n}+\big\|A_{m}-A_{m}^{(k)}\big\|\,\|x_{m}\|_{E_n} +\big\|A_{n}^{(k)}x_{n}-A_{m}^{(k)}x_{m}\big\|_{E}\\
 &< \frac{\eps}{3M}M+\frac{\eps}{3M}M+\frac{1}{k}\leq \eps.
\qedhere
 \end{align*}
\end{proof}

Now we assume that the underlying spaces are Hilbert spaces.
It is well-known that an operator is compact if and only if so is its adjoint (see e.g.\ \cite[Theorem~III.4.10]{kato}).
We prove an analogous result for discrete compactness of a sequence of operators.
\begin{prop}\label{propdisccompofadj}
 Let 
$A\in L(H)$ and $A_n\in L(H_n)$, $n\in\N$.
If
$$A_nP_n\slong AP, \quad A_n^*P_n \slong A^*P ,\quad n\to\infty,$$ 
then the following are equivalent:
\begin{enumerate}[label=\rm{\roman{*})}]
 \item the sequence $\left(A_n\right)_{n\in\N}$ is discretely compact;
 \item the sequence $\left(A_n^*A_n\right)_{n\in\N}$ is discretely compact;
 \item the sequence $\left(A_n^*\right)_{n\in\N}$ is discretely compact.
\end{enumerate}
\end{prop}

\begin{proof}
The claims ``i)~$\To~$ii)'' and ``iii)~$\To$~ii)'' follow from Lemma~\ref{lemmapertdisccomp}~i),~ii).
It is left to show the claim ``ii)~$\To$~i)''; the proof of ``ii)~$\To$~iii)'' is analogous.

Let $M>0$ and let $I\subset\N$ be an infinite subset such that there exist $y_n\in H_n$, \linebreak
$n\in I$, with $\|y_n\|_{H_n}\leq M$, $n\in I$.
We show that there exists a convergent subsequence of $(A_ny_n)_{n\in I}\subset H_0$ with limit in $H$.

Since $H_0$ is weakly compact, there exists an infinite subset $I_2\subset I$ such that $(y_n)_{n\in I_2}\subset H_0$ is weakly convergent; denote the weak limit by $y$.
Since $P_n\s P$, it is easy to see that $y=Py\in H$.
In addition, $A_n^*P_n\s A^*P$ implies $$A_ny_n\stackrel{w}{\tolong}Ay, \quad n\in I_2, \quad n\to\infty.$$
Below we show that $(\|A_ny_n\|_{H_0})_{n\in I_3}$ converges to $\|Ay\|_{H_0}$ for some infinite subset $I_3\subset I_2$; then we obtain the desired convergence $\|A_ny_n-Ay\|_{H_0}\to 0$ as $n\in I_3$, $n\to\infty$.

By the assumptions, the sequence $\left(A_n^*A_n\right)_{n\in\N}$ is discretely compact, thus there exist an infinite subset $I_3\subset I_2$ and $x\in H$ such that $\left(A_n^*A_ny_n\right)_{n\in I_3}$ converges in $H_0$ to~$x$.
On the other hand, the strong convergences $A_nP_n\s AP$ and $A_n^*P_n\s A^*P$ imply the weak convergence 
$$A_n^*A_ny_n\stackrel{w}{\tolong}A^*Ay, \quad n\in I_3, \quad n\to\infty.$$ 
By the uniqueness of the weak limit, we have $A^*Ay=x$.
So we obtain, for $n\in I_3$,
$$\|A_ny_n\|_{H_0}^2=\langle A_n^*A_ny_n,y_n\rangle_{H_0}\tolong\langle x,y\rangle_{H_0}=\langle A^*Ay,y\rangle_{H_0}=\|Ay\|_{H_0}^2,\quad n\to\infty;$$
this proves the claim.
\end{proof}

Related to discrete compactness is the notion of {collectively compact} sets of bounded linear operators (see  Anselone and Palmer \cite{anselonepalmer}).

\begin{definition}\label{defcollectivelycomp}
 Let $\mathcal{B}$ be  the closed unit ball in $H_0$. 
 A subset $\mathcal{K}\subset L(H_0)$ is called \emph{collectively compact}\index{compactness!collectively compact operator set} if the set 
$$\mathcal K\mathcal B=\left\{Kx:\,K\in\mathcal K, \,x\in \mathcal B\right\}\subset H_0$$ is relatively compact in $H_0$.
\end{definition}

\begin{rem}\label{remcollectivelycompact}
  \begin{enumerate}[label=\rm{\roman{*})}]
   \item Every operator of a collectively compact set is compact.
   \item A set $\{A_n:\,n\in\N\}$ is collectively compact if and only if the operators $A_n,\,n\in\N$, are compact and the sequence $(A_n)_{n\in\N}$ is discretely compact. 
 \end{enumerate}
 \end{rem}

The following result 
yields sufficient conditions on Hilbert space operators such that strong convergence implies norm convergence.

\begin{prop}\label{corgsrimpliesgnr}
 Let 
$A\in L(H)$ and $A_n\in L(H_n)$, $n\in\N$. 
Assume that $A_n$, $n\in\N,$ are all compact operators and $(A_n)_{n\in\N}$ is a discretely compact sequence.
If $A_nP_n\s AP$ and $A_n^*P_n\s A^*P$, then $A$ is compact and $A_nP_n\to AP$.
\end{prop}
\begin{proof}
 It is well-known that if $A_n,\,n\in\N,$ are all compact operators, then so are $A_nP_n$, $A_n^*P_n,\,n\in\N$.
Proposition~\ref{propdisccompofadj} yields the discrete compactness of the sequence $(A_n^*)_{n\in\N}$.
By Lemma~\ref{lemmapertdisccomp}~i), the sequences  $(A_nP_n)_{n\in\N}$,  $(A_n^*P_n)_{n\in\N}$ are discretely compact.
Remark~\ref{remcollectivelycompact}~ii) implies that $\{A_nP_n:\,n\in\N\}$, $\{A_n^*P_n:\,n\in\N\}$ are collectively compact sets.
Then, by \cite[Proposition 2.1 (a)$\,\, \To\,~$(b)]{anselonepalmer}, so are the sets $\{A_nP_n-AP:\,n\in\N\}$, $\{A_n^*P_n-A^*P:\,n\in\N\}$, and $AP$ (and thus $A$) is a compact operator.
Now the claim follows from \cite[Theorem 3.4~(c)]{anselonepalmer}.
\end{proof}

\subsection{Proofs of main results}\label{subsectionproofs}
In this subsection we prove the theorems in Subsection~\ref{subsectionmainresults}.

The following two elementary results will be used later on.
\begin{lemma}\label{AndAJndJnew}
Let $T\in C(E)$ and $T_n\in C(E_n)$, $n\in\N$. 
Assume that $T_n\gsr T$. Then
for all $x\in \dom(T)$ there exists a sequence of elements $x_n\in\dom(T_n), \,n\in\N,$ such that 
    \beq \|x_n-x\|_{E_0}+\|T_nx_n-Tx\|_{E_0}\tolong 0, \quad n\to\infty.\label{eq.disc.conv}\eeq
\end{lemma}

\begin{proof}
By Definition~\ref{defresconv}~i) of $T_n\gsr T$, there exist $n_0\in\N$ and $\lm\in\underset{n\geq n_0}{\bigcap}\rho(T_n)\cap\rho(T)$ such that
\beq (T_n-\lm)^{-1}P_n\slong (T-\lm)^{-1}P, \quad n\to\infty.\label{resconv}\eeq
Let $x\in\dom(T)$ and define $$x_n:=(T_n-\lm)^{-1}P_n(T-\lm)x\in\dom(T_n), \quad n\geq n_0.$$
Then, using $P_n\s P$ and \eqref{resconv}, it is easy to verify that~\eqref{eq.disc.conv} holds.
\end{proof}

\begin{lemma}\label{lemma.regofbdd.compact}
Let $T_n\in C(E_n)$, $n\in\N$, and let $K\subset\Delta_b\big((T_n)_{n\in\N}\big)$ be a compact subset.
Then there exist $M_{K}>0$ and $n_{K}\in\N$ such that
$$\forall\,\lm\in K:\quad \lm\in\rho(T_n), \quad \|(T_n-\lm)^{-1}\|\leq M_K, \quad n\geq n_K.$$
\end{lemma}

\begin{proof}
Assume that the claim is false, i.e.\ no such $M_K>0$ exists. 
Then there exist an infinite subset $I_1\subset\N$ and $(\lm_n)_{n\in I_1}\subset K$ such that $\|(T_n-\lm_n)^{-1}\|\to\infty$ as $n\in I_1$, $n\to\infty$.
By the compactness of $K$, there are $\lm\in K$ and an infinite subset $I_2\subset I_1$ so that $(\lm_n)_{n\in I_2}$ converges to $\lm$.
Since $\lm\in \Delta_b\big((T_n)_{n\in\N}\big)$, there exist $M_{\lm}>0$ and an infinite subset $I_3\subset I_2$ such that $\lm\in\rho(T_n)$ and $\|(T_n-\lm)^{-1}\|\leq M_{\lm}$ for all $n\in I_3$.
Then, for every $n\in I_3$ so large that $|\lm_n-\lm|\leq 1/(2M_{\lm})$, a Neumann series argument yields
\begin{align*}
\|(T_n-\lm_n)^{-1}\|&=\big\|(T_n-\lm)^{-1}\big(I-(\lm_n-\lm)(T_n-\lm)^{-1}\big)^{-1}\big\|\leq 2 M_{\lm}.
\end{align*}
The obtained contradiction proves the claim.
\end{proof}

For generalized strong/norm resolvent convergence we assume the resolvents to converge for one particular $\lm$.
In the following result we investigate for which points the resolvents then converge as well.
\begin{prop}\label{propgsrgnrregionofbdd}
Let $T\in C(E)$ and $T_n\in C(E_n)$, $n\in\N$. 
\begin{enumerate}[label=\rm{\roman{*})}]
\item Assume that $T_n\gsr T$. Then
$\Delta_b\big((T_n)_{n\in\N}\big)\subset\C\backslash \sigma_{\rm app}(T)$ and, for any $\lm\in\Delta_b\big((T_n)_{n\in\N}\big)\cap\rho(T)$, \beq (T_n-\lm)^{-1}P_n\slong(T-\lm)^{-1}P, \quad n\to\infty.\label{eqconvtoshow}\eeq
\item Assume that $T_n\gnr T$. Then $\rho(T)\subset\Delta_b\big((T_n)_{n\in\N}\big)$ and, for any $\lm\in\rho(T)$, \beq (T_n-\lm)^{-1}P_n\tolong(T-\lm)^{-1}P, \quad n\to\infty.\label{eqconvtoshow2}\eeq
\end{enumerate}
\end{prop}

\begin{proof}
i) By Definition~\ref{defresconv}~i) of $T_n\gsr T$, there exist $n_0\in\N$ and an element $\lm_0\in\underset{n\geq n_0}{\bigcap}\rho(T_n)\cap\rho(T)$ such that \vspace{-2mm}
\beq (T_n-\lm_0)^{-1}P_n\slong (T-\lm_0)^{-1}P, \quad n\to\infty.\label{eqRnslongR}\eeq
Let $\lm\in\Delta_b\big((T_n)_{n\in\N}\big)$. 
By Definition~\ref{defresconv}~iii) of the region of boundedness, there exists $n_1\in\N$ (without loss of generality $n_1\geq n_0$) 
such that $\lm\in\rho(T_n)$, $n\geq n_1$, and $M:=\sup_{n\geq n_1}\|(T_n-\lm)^{-1}\|<\infty$.
First assume that $\lm$ belongs to the approximate point spectrum $\sigma_{\rm app}(T)$. Then there exists $x\in\dom(T)$ with $\|x\|_E=1$ and $\|(T-\lm)x\|_E<1/(2M)$.
By Lemma~\ref{AndAJndJnew}, there exists a sequence of elements $x_n\in\dom(T_n)$, $n\in\N$, such that $\|x_n\|_{E_n}=1$ and $\|(T_n-\lm)x_n\|_{E_n}<1/(2M)$ for all large enough $n\in\N$.
The obtained contradiction to $M=\sup_{n\geq n_1}\|(T_n-\lm)^{-1}\|$ implies $\lm\notin\sigma_{\rm app}(T)$.

Now assume that $\lm\in\Delta_b\big((T_n)_{n\in\N}\big)\cap\rho(T)$.
If we set, for $n\geq n_1$,
\begin{align*} S_n(\lm):=&\,I+(\lm_0-\lm)(T_n-\lm_0)^{-1}P_n =\,(T_n-\lm)(T_n-\lm_0)^{-1}P_n+(I-P_n),\end{align*}
then a straightforward application of the first resolvent identity yields
\beq \begin{aligned}
&S_n(\lm)\big((T_n-\lm)^{-1}P_n-(T-\lm)^{-1}P\big)\\
&=\big((T_n-\lm_0)^{-1}P_n-(T-\lm_0)^{-1}P\big)\big(I+(\lm-\lm_0)(T-\lm)^{-1}P\big), \quad n\geq n_1. 
\end{aligned}\label{eqSnRnDiff}
\eeq
The operators $S_n(\lm)$, $n\geq n_1$, are invertible,
\beq \begin{aligned} \hspace{-1.5mm}S_n(\lm)^{-1}&=(T_n-\lm_0)(T_n-\lm)^{-1}P_n+(I-P_n)=I+(\lm-\lm_0)(T_n-\lm)^{-1}P_n,\hspace{-2mm}\end{aligned}\label{eqSninverse}\eeq
and the inverses are uniformly bounded  since $\lm\in\Delta_b\big((T_n)_{n\in\N}\big)$.
Now the claimed convergence~\eqref{eqconvtoshow} follows from~\eqref{eqRnslongR} and~\eqref{eqSnRnDiff}.

ii)
Let $\lm\in\rho(T)$.
By Definition~\ref{defresconv}~ii) of $T_n\gnr T$, there exist $n_0\in\N$  and an element 
$\lm_0\in\underset{n\geq n_0}{\bigcap}\rho(T_n)\cap\rho(T)$ such that \vspace{-2mm}
\beq (T_n-\lm_0)^{-1}P_n\tolong (T-\lm_0)^{-1}P, \quad n\to\infty.\label{eqRntolongR}\eeq
This implies 
\begin{align*} S_n(\lm)=&\, I+(\lm_0-\lm)(T_n-\lm_0)^{-1}P_n \tolong I+(\lm_0-\lm)(T-\lm_0)^{-1}P=:S, \quad n\to\infty.\end{align*}
Since $S$ is boundedly invertible, \cite[Theorem~IV.1.16]{kato} yields the existence of some $n_1\in\N$ 
such that the operators $S_n(\lm)$, $n\geq n_1$, are uniformly boundedly invertible.
Then $(T_n-\lm)$, $n\geq n_1$, are uniformly boundedly invertible because, by~\eqref{eqSninverse}, $$(T_n-\lm)^{-1}=(\lm-\lm_0)^{-1}\left.\big(S_n(\lm)^{-1}-I\big)\right|_{E_n}, \quad n\geq n_1.$$
Now, analogously as in~i), the claim~\eqref{eqconvtoshow2} follows from~\eqref{eqRntolongR} and~\eqref{eqSnRnDiff}.
\end{proof}

Now we are ready to prove the main results of Subsection~\ref{subsectionmainresults}.

\begin{proof}[Proof of Theorem~{\rm\ref{thmlocal}}]
i) 
Let $\lm\in\sigma(T)$ and $\eps>0$ satisfy~\eqref{eqneighbourhoodinregionofbdd}.
Choose $\delta>0$ with $\delta<\eps$. Assume that there exists an infinite subset $I\subset\N$ with ${\rm dist}(\lm,\sigma(T_n))\geq\delta$, $n\in I$. 
For $\Gamma:=\partial B_{{\delta}/2}(\lm)$ define the contour integrals
\begin{align*}
 P_{\Gamma}:=\frac{1}{2\pi\I}\int_{\Gamma} (T-z)^{-1}\,\rd z, \quad  P_{\Gamma,n}:=\frac{1}{2\pi\I}\int_{\Gamma} (T_n-z)^{-1}\,\rd z, \quad n\in I.
\end{align*}
The operator $P_{\Gamma}$ is the spectral projection corresponding to $\lm\in\sigma(T)$. However, since $z\mapsto (T_n-z)^{-1}$ is holomorphic in $B_{\delta}(\lm)$, we have $P_{\Gamma,n}=0$, $n\in I$.
Let $x\in E_0$ be arbitrary. For $n\in I$ define the function $f_n:\Gamma\to [0,\infty)$ by $f_n(z):=\|(T-z)^{-1}Px-(T_n-z)^{-1}P_nx\|_{E_0}$.
Then $$\|P_{\Gamma}Px-P_{\Gamma,n}P_nx\|_{E_0}\leq \frac{1}{2\pi}\int_{\Gamma} f_n (z)\,\rd |z|, \quad n\in I.$$
Note that $f_n(z)\to 0$, $n\to\infty$, for every $z\in\Gamma$, and $f_n,\,n\in\N,$ are uniformly bounded by the compactness of $\Gamma\subset\Delta_b\big((T_n)_{n\in\N}\big)$ and by Lemma~\ref{lemma.regofbdd.compact}.
Lebesgue's dominated convergence theorem implies $\|P_{\Gamma}Px-P_{\Gamma,n}P_nx\|_{E_0}\to 0$ as $n\to\infty$.
Hence $P_{\Gamma,n}P_n \s P_{\Gamma}P, \,n\to\infty,$ and so we arrive at the contradiction $P_{\Gamma}=0$. 
Therefore, there exists $n_{\delta}\in\N$ such that ${\rm dist}(\lm,\sigma(T_n))<\delta$, $n\geq n_{\delta}$.
Since $\delta$ can be chosen arbitrarily small, we finally obtain ${\rm dist}(\lm,\sigma(T_n))\to 0$, $n\to\infty$.

ii)
Let $\lm\in\rho(T)\cap\Delta_b\left((T_n)_{n\in\N}\right)$. 
Definition~\ref{defresconv}~iii) of the region of boundedness implies that there exist $n_0\in\N$ and $M>0$ such that $\lm\in\rho(T_n),\,n\geq n_0$, 
and $\left\|(T_n-\lm)^{-1}\right\|\leq M$, $n\geq n_0$.
As a consequence, $$\dist(\lm,\sigma(T_n))\geq \frac{1}{\left\|(T_n-\lm)^{-1}\right\|}\geq \frac{1}{M}, \quad n\geq n_0,$$
so $\lm$ cannot be the limit of a sequence of points in the spectra of $T_n,\,n\geq n_0$.
\end{proof}

\begin{proof}[Proof of Theorem~{\rm\ref{thmgeneralizedresconv}}]
i)
By Proposition~\ref{propgsrgnrregionofbdd}~ii), we have $\rho(T)\subset\Delta_b\big((T_n)_{n\in\N}\big)$. Now the claim follows from Theorem~\ref{thmlocal}~ii).

ii)
Since $T$ is assumed to be selfadjoint, it satisfies $\sigma(T)=\sigma_{\rm app}(T)$ and thus Proposition~\ref{propgsrgnrregionofbdd}~i) implies $\Delta_b\big((T_n)_{n\in\N}\big)\subset\rho(T)$.
In addition, since $T_n$ is selfadjoint, we have ${\rm dist}(\lm,\sigma(T_n))=\|(T_n-\lm)^{-1}\|^{-1}$ for any $\lm\in\rho(T_n)$, which implies~\eqref{eq.Deltab.selfadjoint}. 
\end{proof}

\begin{proof}[Proof of Theorem~{\rm\ref{mainthmforgsr}}]
i) Since $T$ has compact resolvent, it satisfies $\sigma(T)=\sigma_p(T)=\sigma_{\rm app}(T)$.
Proposition~\ref{propgsrgnrregionofbdd}~i)  implies that $\Delta_b\big((T_n)_{n\in\N}\big)\subset\rho(T)$.

Conversely, take $\lm\in \C\backslash\Delta_b\big((T_n)_{n\in\N}\big)$. Note that $\sigma(T_n)=\sigma_p(T_n)$ since $T_n$ is assumed to have compact resolvent.
Then there are an infinite subset $I\subset \N$ and $x_n\in \dom(T_n)$, $n\in I$, with 
\begin{equation}\label{eq.reg.seq}
\|x_n\|_{E_n}=1, \quad n\in I, \quad \|(T_n-\lm)x_n\|_{E_n}\tolong 0, \quad n\to\infty.
\end{equation}
Define $y_n:=(T_n-\lm_0)x_n$ for $n\in I.$
Then $(\|y_n\|_{E_n})_{n\in I}$ is a bounded sequence.
Since  $\left((T_n-\lm_0)^{-1}\right)_{n\in\N}$ is discretely compact by the assumptions, there exist
$x\in E$ and an infinite subset $\widetilde I\subset I$ so that $\|x_n-x\|_{E_0}\to 0$ as $n\in\widetilde I$, $n\to\infty$.
By~\eqref{eq.reg.seq}, we have $\|x\|_E=1$ and $\|y_n-(\lm-\lm_0)x\|_{E_0}\to 0$.
However, $(T_n-\lm_0)^{-1}P_n\s (T-\lm_0)^{-1}P$ then yields
$$x_n=(T_n-\lm_0)^{-1}y_n\tolong (\lm-\lm_0)(T-\lm_0)^{-1}x\in\dom(T),\quad n\in \widetilde I,\quad  n\to\infty.$$ 
By the uniqueness of the limit, we obtain $x\in\dom(T)$ and $Tx=\lm x$.
Since $x\neq 0$, we have $\lm\in\sigma(T)$.

The convergence~\eqref{eq.gsr.atalllm} for all $\lm\in\rho(T)$ now follows from Proposition~\ref{propgsrgnrregionofbdd}~i).

ii)
Spectral exactness follows from claim~i) and Theorem~\ref{thmlocal}; note that all $\lm\in\sigma(T)$ are isolated since $T$ is assumed to have compact resolvent.
 In an analogous way as in the proof of~\ref{thmlocal}~i), one may prove that the corresponding spectral projections converge strongly, $P_{\Gamma,n}P_n\s P_{\Gamma}P$. This implies that for $x=P_{\Gamma}Px$ in the algebraic eigenspace of $\lm$ there exists a sequence of elements $x_n:=P_{\Gamma,n}P_nx\in\ran(P_{\Gamma,n})$, $n\in\N$, with $\|x_n-x\|_{E_0}\to 0$, and the normalized elements converge as well. This proves that $m={\rm rank}\,P_{\Gamma}\leq\liminf_{n\to\infty}{\rm rank}\,P_{\Gamma,n}$. 

To prove that $\limsup_{n\to\infty}{\rm rank}\,P_{\Gamma,n}\leq m$ (and thus ${\rm rank}\,P_{\Gamma,n}=m$ for all sufficiently large $n$), let $\lm_n\in\sigma(T_n)$, $n\in\N$, such that $\lm_n\to\lm\in\sigma(T)$ as $n\to\infty$.
For $n\in\N$ denote by  $k_n$ the ascent of $\lm_n$, i.e.\ the smallest $k\in\N$ such that $(T_n-\lm_n)^kP_{\Gamma,n}=(T_n-\lm_n)^{k+1}P_{\Gamma,n}$;
then there exist $k_n$ orthonormal elements $x_n^{(k)}\in\ran(P_{\Gamma,n})$ with $(T_n-\lm_n)x_n^{(1)}=0$ and $(T_n-\lm_n)x_n^{(k)}=x_n^{(k-1)}$ for $k=2,\dots,k_n$.
By induction over $k\in\N$, we prove that if there exists an infinite subset $I^{(k)}\subset\N$ such that $k\leq k_n$, $n\in I^{(k)}$, then there exist $x^{(k)}\in\ran(P_{\Gamma})$ and an infinite subset $\widetilde I^{(k)}\subset I^{(k)}$ such that $\|x^{(k)}\|_E=1$, $(T-\lm)^kx^{(k)}=0$ and  $\|x_n^{(k)}-x^{(k)}\|_{E_0}\to 0$  as $n\in \widetilde I^{(k)}$, $n\to\infty$; then we obtain $\limsup_{n\to\infty}{\rm rank}\,P_{\Gamma,n}\leq m$.

If $k=1$, then the claim follows from the proof of  i) since $\lm_n\to\lm$ and thus $x_n^{(1)}$, $n\in I^{(1)}$, satisfy~\eqref{eq.reg.seq}.
If $k>1$, set $y_n^{(k)}:=(T_n-\lm_0)x_n^{(k)}=x_n^{(k-1)}+(\lm_n-\lm_0)x_n^{(k)}$, $n\in I^{(k)}$. 
Since $\lm_n\to\lm$, the sequence $(\|y_n^{(k)}\|_{E_n})_{n\in I^{(k)}}$ is bounded. By proceeding as in i) and using the induction hypothesis, we obtain the claim.
\end{proof}

\begin{rem}
Theorem~\ref{mainthmforgsr} would also follow from Stummel's results \cite[S\"atze 2.2.(1), 3.2.(8)]{stummel2} applied to the bounded operators $A,B\in L(D,E)$ and $A_n,B_n\in L(D_n,E_n)$, $n\in\N$, where the Banach spaces
$D:=\dom(T)$, $D_n:=\dom(T_n)$
are equipped with the graph norm of $T$, $T_n$, respectively,
 the operators $B:D\to E, \,B_n:D_n\to E_n$ are the natural embeddings, and $A,A_n$ are defined by
$Ax:=Tx$, $A_nx_n:=T_nx_n.$
However, it is very technical to check that the assumptions of Stummel's results are satisfied (see~\cite[Chapter~1]{boegli-phd} for this approach).
\end{rem}

\begin{proof}[Proof of Theorem~{\rm\ref{thmgsrimpliesgnrbasic}}]
The claim for $\lm=\lm_0$ is an immediate consequence of Proposition~\ref{corgsrimpliesgnr} applied to
$A=(T-\lm_0)^{-1}$, $A_n=(T_n-\lm_0)^{-1}$, $n\in\N.$
The generalized norm resolvent convergence for every $\lm\in\rho(T)$ follows from  Proposition~\ref{propgsrgnrregionofbdd}~ii).
\end{proof}

\section{Generalized strong resolvent (gsr) convergence}\label{sectionsuffcondforstrongresconv}
In this section we establish sufficient conditions for gsr-convergence of a sequence of linear operators in varying Banach spaces.
First, we find conditions that directly yield gsr-convergence,
but then we also give sufficient conditions for gsr-convergence of a sequence of operators $A_n=T_n+S_n\in C(E_n),\,n\in\N,$ if we know it for the operators $T_n,\,n\in\N$ (see Subsection~\ref{subsectiondirectgsr}).
Afterwards, in Subsection~\ref{subsection2x2matricesconv}, we derive sufficient conditions on a sequence of $2\times 2$ block operator matrices that imply gsr-convergence.
Then we consider gsr-convergence of infinite operator matrices (see Subsections~\ref{subsectionconvlemmas},~\ref{subsectioninfititegsr}).

\subsection{Direct criteria and perturbation results}\label{subsectiondirectgsr}
Let $E_0$ be a Banach space and $E, \,E_n\subset E_0, \,n\in\N,$ be closed subspaces. Denote by 
$P:E_0\to E$, $P_n:E_0\to E_n$, $n\in\N,$ projections onto the respective subspaces that converge strongly, $P_n\s P$.

The next proposition is 
a generalization of Weidmann's result~\cite[Satz~9.29~a)]{weid1} who considers selfadjoint operators.
\begin{theorem} \label{propcoreforgsr}
Let $T\in C(E)$ and  $T_n\in C(E_n), \,n\in\N$.
 Let $\Phi\subset\dom(T)$ be a core of $T$ such that 
for all $x\in\Phi$ there exists $n_0(x)\in\N$ with 
$$\forall\,n\geq n_0(x):\quad P_nx\in\dom(T_n), \quad \|T_nP_nx- Tx\|_{E_0}, \quad n\to\infty.$$ 
Suppose that $\Delta_b\left((T_n)_{n\in\N}\right)\cap\rho(T)\neq\emptyset$. Then $T_n\gsr T$.
\end{theorem}

\begin{proof}
The proof of Weidmann's result~\cite[Satz~9.29~a)]{weid1} remains valid in the non-selfadjoint case.
The only place in the latter result where the selfadjointness of $T\in C(E)$, $T_n\in C(E_n), \,n\in\N$, is used is the consequence 
\beq \Delta_b\left((T_n)_{n\in\N}\right)\cap\rho(T)\neq\emptyset \label{eqintersectionnonempty}\eeq
(since $\C\backslash\R$ belongs to the intersection).
Since here $T\in C(E)$, $T_n\in C(E_n), \,n\in\N$, are not assumed to be selfadjoint, we require~\eqref{eqintersectionnonempty} to be satisfied by assumption.
\end{proof}

For bounded operators strong convergence implies gsr-convergence.
\begin{lemma}\label{lemmabddstongconv}
 Let $B\in L(E)$ and $B_n\in L(E_n)$, $n\in\N$, satisfy $B_nP_n\s BP$. Then \beq \Big\{\lm\in\C:\,|\lm|>\limsup_{n\to\infty}\|B_n\|\Big\}\subset\Delta_b\big((B_n)_{n\in\N}\big),\vspace{-1mm}\label{eqDeltabforbdd}\eeq
and $$\forall\,\lm\in \Delta_b\big((B_n)_{n\in\N}\big)\cap\rho(B):\quad (B_n-\lm)^{-1}P_n\slong (B-\lm)^{-1}P, \quad n\to\infty.$$
\end{lemma}

\begin{proof}
Let $\lm\in\C$ satisfy $|\lm|>\limsup_{n\to\infty}\|B_n\|=:M$. Then there exist $\eps>0$ and $n_0\in\N$ with $$|\lm|>\|B_n\|+\eps, \quad n\geq n_0.$$
Now a Neumann series argument yields that, for $n\geq n_0$, the operator $(B_n-\lm)$ is invertible, with 
$$(B_n-\lm)^{-1}=-\lm^{-1}(I-\lm^{-1}B_n)^{-1}, \quad \|(B_n-\lm)^{-1}\|\leq \frac{|\lm^{-1}|}{1-|\lm^{-1}|\|B_n\|}\leq \frac{1}{\eps}.$$
This proves the inclusion~\eqref{eqDeltabforbdd}.

 Let $\lm\in \Delta_b\big((B_n)_{n\in\N}\big)\cap\rho(B)$.
Then there exists $n_0\in\N$ such that $\lm\in\rho(B_n)$, $n\geq n_0$, and $\big((B_n-\lm)^{-1}\big)_{n\geq n_0}$ is a bounded sequence.
Let $y\in E_0$ and define $x:=(B-\lm)^{-1}Py$. It is easy to verify that, for every $n\geq n_0$,
\begin{align}
 &\big((B_n-\lm)^{-1}P_n-(B-\lm)^{-1}P\big)y\notag\\
&=(B_n-\lm)^{-1}P_n(BP-B_nP_n)x-(B_n-\lm)^{-1}P_n (P-P_n)y-(P-P_n)x.\label{eqBnBresolvents}
\end{align}
Since the sequence $\big(\|(B_n-\lm)^{-1}\|\big)_{n\geq n_0}$ is bounded,
the assumptions $B_nP_n\s BP$ and $P_n\s P$ imply that the right-hand side of~\eqref{eqBnBresolvents} converges to $0$. 
\end{proof}

Now we consider sums $A=T+S$ and $A_n=T_n+S_n$, $n\in\N$.
We study perturbation results for generalized strong resolvent convergence, i.e.\ we establish sufficient conditions that $T_n\gsr T$ implies $A_n\gsr A$.

\begin{theorem}\label{perturbgsr}
Let $T\in C(E)$ and $T_n\in C(E_n), \,n\in\N$.
Let  $S$ and $S_n$, $n\in\N$, be linear operators in $E$ and $E_n$, $n\in\N$, with $\dom(T)\subset\dom(S)$ and $\dom(T_n)\subset\dom(S_n)$, $n\in\N$, respectively.
Define $$A:=T+S, \quad A_n:=T_n+S_n, \quad n\in\N.$$
Suppose that there exist $\lm\in\underset{n\in\N}{\bigcap}\rho(T_n)\cap\rho(T)$ 
and $\gamma_{\lm}<1$ with
\begin{align}
\left\|S(T-\lm)^{-1}\right\|<1,\quad \left\|S_n(T_n-\lm)^{-1}\right\|\leq \gamma_{\lm},\quad n\in\N.\label{eqSTbdd}
\end{align}
If
\beq\begin{array}{rl}(T_n-\lm)^{-1}P_n \hspace{-2mm}&\slong\, (T-\lm)^{-1}P,\\[1mm] S_n(T_n-\lm)^{-1}P_n \hspace{-2mm}&\slong\, S(T-\lm)^{-1}P,\end{array}\quad n\to\infty,\label{eqconvgsrT}\eeq
then
$\lm\in\underset{n\in\N}{\bigcap}\rho(A_n)\cap\rho(A)$ and  $(A_n-\lm)^{-1}P_n\s (A-\lm)^{-1}P$ as $n\to\infty.$
\end{theorem}

\begin{rem}
 The inequalities~\eqref{eqSTbdd} imply that $S$ is $T$-bounded with $T$-bound $<1$ and $S_n$ is $T_n$-bounded with $T_n$-bound $\leq\gamma_{\lm}<1$. 
\end{rem}

\begin{proof}[Proof of Theorem~{\rm\ref{perturbgsr}}]
 Let $\lm$ and $\gamma_{\lm}<1$ satisfy the assumptions.
For $n\in\N$, by a Neumann series argument, $(A_n-\lm)$ is boundedly invertible and
\beq
\begin{aligned}\label{neumanngamma}
(A_n-\lm)^{-1}&=(T_n-\lm)^{-1}\left(I+S_n(T_n-\lm)^{-1}\right)^{-1}, \\ 
\big\|\left(I+S_n(T_n-\lm)^{-1}\right)^{-1}\big\|&\leq\frac{1}{1-\gamma_{\lm}}.
\end{aligned}
\eeq
Analogously we obtain 
\beq (A-\lm)^{-1}=(T-\lm)^{-1}\left(I+S(T-\lm)^{-1}\right)^{-1}.\label{neumann1}\eeq
Thus $\lm\in\underset{n\in\N}{\bigcap}\rho(A_n)\cap\rho(A)$. Let $x\in E_0$ and define
$y:=\left(I+S(T-\lm)^{-1}\right)^{-1}Px\in E$. 
Since $\big((T_n-\lm)^{-1}\big)_{n\in\N}$ is strongly convergent, $C:=\sup_{n\in\N}\|(T_n-\lm)^{-1}\|<\infty$.
Then 
\begin{align*}
 &\left\|(A_n-\lm)^{-1}P_nx-(A-\lm)^{-1}Px\right\|_{E_0}\\
&\leq \left\|(T_n-\lm)^{-1}P_n\left(\left(I+S_n(T_n-\lm)^{-1}\right)^{-1}P_nx-\left(I+S(T-\lm)^{-1}\right)^{-1}Px\right)\right\|_{E_0}\\
&\quad +\left\|\big((T_n-\lm)^{-1}P_n-(T-\lm)^{-1}P\big)\left(I+S(T-\lm)^{-1}\right)^{-1}Px\right\|_{E_0}\\
&\leq C \left\|\left(I+S_n(T_n-\lm)^{-1}\right)^{-1}P_nx-\left(I+S(T-\lm)^{-1}\right)^{-1}Px\right\|_{E_0}\\
&\quad +\left\|\big((T_n-\lm)^{-1}P_n-(T-\lm)^{-1}P\big)y\right\|_{E_0}.
\end{align*}
The first convergence in~\eqref{eqconvgsrT} yields $$\left\|\big((T_n-\lm)^{-1}P_n-(T-\lm)^{-1}P\big)y\right\|_{E_0}\tolong 0, \quad n\to\infty.$$
By \eqref{neumanngamma} and \eqref{neumann1}, we have $-1\in \Delta_b\left((S_n(T_n-\lm)^{-1})_{n\in\N}\right)\cap\rho(S(T-\lm)^{-1})$.
Hence, the second convergence in~\eqref{eqconvgsrT} implies that, by Lemma~\ref{lemmabddstongconv},
$$\left\|\left(I+S_n(T_n-\lm)^{-1}\right)^{-1}P_nx-\left(I+S(T-\lm)^{-1}\right)^{-1}Px\right\|_{E_0}\tolong 0,\quad n\to\infty.$$
Altogether, we have $(A_n-\lm)^{-1}P_n\s(A-\lm)^{-1}P, \, n\to\infty.$
 \end{proof}

The following result is an immediate consequence of Theorem \ref{perturbgsr} for the case that the perturbations are bounded operators.
\begin{corollary} \label{corperturbgsr}
Let $T\in C(E)$ and $T_n\in C(E_n), \,n\in\N$.
Let  $S\in L(E)$ and $S_n\in L(E_n)$, $n\in\N$.
Define $$A:=T+S, \quad A_n:=T_n+S_n, \quad n\in\N.$$
Suppose that, for some $\lm\in\underset{n\in\N}{\bigcap}\rho(T_n)\cap\rho(T)$ and some $\gamma_{\lm}<1$, we have 
\begin{align*}\label{eqSSnrelbdd}
\left\|S\right\|<\|(T-\lm)^{-1}\|^{-1},\quad 
\left\|S_n\right\|\leq \gamma_{\lm}\,\|(T_n-\lm)^{-1}\|^{-1},\quad n\in\N,
\end{align*}
and 
\begin{align*}
 (T_n-\lm)^{-1}P_n\slong (T-\lm)^{-1}P, \quad S_nP_n\slong SP,\quad &n\to\infty.
\end{align*}
Then
$\lm\in\underset{n\in\N}{\bigcap}\rho(A_n)\cap\rho(A)$ and  $(A_n-\lm)^{-1}P_n\s (A-\lm)^{-1}P$ as $n\to\infty.$
\end{corollary}

\subsection{Results for $2\times 2$ block operator matrices}\label{subsection2x2matricesconv}
In this subsection we consider diagonally dominant $2\times 2$ operator matrices. To this end, all spaces are assumed to consist of two components.
For $i=1,2$, let $E_i^{(0)}$ be a Banach space and let $E_i, \,E_i^{(n)}\subset E_i^{(0)}, \,n\in\N,$ be closed subspaces;
denote by $P_i:E_i^{(0)}\to E_i$, $P_i^{(n)}:E_i^{(0)}\to E_i^{(n)}$, $n\in\N,$ projections onto the respective subspaces that converge strongly, $P_i^{(n)}\s P_i$, $n\to\infty$.
In the product space $\cE:=E_1\oplus E_2$ we consider a $2\times 2$ block operator matrix\index{block operator matrix!$2\times 2$}
\be\cA:=\bmat A & B\\C & D \emat, \quad \dom(\cA):=\left(\dom(A)\cap\dom(C)\right)\oplus \left(\dom(B)\cap\dom(D)\right)\ee 
where $A:E_1\to E_1$, $B:E_2\to E_1$, $C: E_1\to E_2$, $D:E_2\to E_2$ are closable operators with dense domains.
We assume that $\cA$ is densely defined.

\begin{definition}{\cite[Definition~2.2.3]{tretter}}\label{defdiagdominant}\index{block operator matrix!diagonally/off-diagonally dominant}
Let $\delta\geq 0$. The block operator matrix $\cA$ is called
  \begin{enumerate}[label=\rm{\roman{*})}] 
   \item \emph{diagonally dominant {\rm(}of order $\delta${\rm)}} if $C$ is $A$-bounded with $A$-bound $\delta_C$, $B$ is $D$-bounded with $D$-bound $\delta_B$, and $\delta=\max\{\delta_C,\delta_B\}$,
   \item \emph{off-diagonally dominant {\rm(}of order $\delta${\rm)}} if $A$ is $C$-bounded with $C$-bound $\delta_A$, $D$ is $B$-bounded with $B$-bound $\delta_D$, and $\delta=\max\{\delta_A,\delta_D\}$.
  \end{enumerate}
\end{definition}

In the product space $\cE^{(n)}:=E_1^{(n)}\oplus E_2^{(n)}$ we consider a $2\times 2$ block operator matrix
\be\cA^{(n)}:=\bmat A^{(n)} & B^{(n)}\\C^{(n)} & D^{(n)} \emat, \quad  \dom(\cA^{(n)}):=\big(\dom(A^{(n)})\cap\dom(C^{(n)})\big)\oplus \big(\dom(B^{(n)})\cap\dom(D^{(n)})\big)\ee 
with operators $A^{(n)}:E_1^{(n)}\to E_1^{(n)}$, $B^{(n)}:E_2^{(n)}\to E_1^{(n)}$, $C^{(n)}: E_1^{(n)}\to E_2^{(n)}$, $D^{(n)}:E_2^{(n)}\to E_2^{(n)}$ that are assumed to be closable with dense domains.
We assume that $\cA^{(n)}$ is densely defined as well.

\begin{theorem}
Suppose that the following holds:
\begin{enumerate}[label=\rm{\roman{*})}]
\item the block operator matrices $\cA$, $\cA^{(n)}$, $n\in\N$, are diagonally dominant;
\item there exists a core $\Phi_1\subset\dom(A)$ of $A$ such that for all $x_1\in\Phi_1$ there exists $n_1(x_1)\in\N$ with the property that
\begin{align*}
&P_1^{(n)}x_1\in\dom(A^{(n)}), \quad n\geq n_1(x_1), \\
&\big\|A^{(n)}P_1^{(n)}x_1- Ax_1\big\|_{E_1^{(0)}}+\big\|C^{(n)}P_1^{(n)}x_1- Cx_1\big\|_{E_1^{(0)}}\tolong 0, \quad n\to\infty;
\end{align*}
\item there exists a core $\Phi_2\subset\dom(D)$ of $D$ such that for all $x_2\in\Phi_2$ there exists $n_2(x_2)\in\N$ with the property that
\begin{align*}
&P_2^{(n)}x_2\in\dom(D^{(n)}), \quad n\geq n_2(x_2), \\
&\big\|D^{(n)}P_2^{(n)}x_2- Dx_2\big\|_{E_2^{(0)}}+\big\|B^{(n)}P_2^{(n)}x_2- Bx_2\big\|_{E_2^{(0)}}\tolong 0, \quad n\to\infty;
\end{align*}
\item  $\Delta_b\left((\cA^{(n)})_{n\in\N}\right)\cap\rho(\cA)\neq\emptyset$.
 \end{enumerate}
Then $\cA^{(n)}\gsr\cA$.
\label{propgsrformatrix}
\end{theorem}

\begin{proof}
First note that $\cA$ is closed by the assumption $\rho(\cA)\neq\emptyset$.
Define the subspace $\Phi:=\Phi_1\oplus\Phi_2\subset\dom(\cA)$. We show that $\Phi$ is a core of $\cA$.

It suffices for each $x:=(x_1,x_2)\in\dom(A)\oplus\dom(D)=\dom(\cA)$ to find a sequence $\big(x^{(m)}\big)_{m\in\N}\subset\Phi$ such that $x^{(m)}\to x$, $\cA x^{(m)}\to\cA x$ as $m\to\infty$.
Since $\Phi_1$ is a core of $A$ and $\Phi_2$ is a core of $D$, there exist sequences
 $\big(x_1^{(m)}\big)_{m\in\N}\subset\Phi_1, \,\big(x_2^{(m)}\big)_{m\in\N}\subset\Phi_2$ with 
\begin{align*}
\begin{array}{ll} 
x_1^{(m)}\tolong x_1, \quad &Ax_1^{(m)}\tolong Ax_1, \\[2mm]
 x_2^{(m)}\tolong x_2, \quad &Dx_2^{(m)}\tolong Dx_2, 
\end{array}
\quad m\to\infty.
\end{align*}
Since $C$ is $A$-bounded and $B$ is $D$-bounded, there exist $a_1, b_1, a_2, b_2\geq 0$ such that
\begin{align*}
\begin{array}{ll} 
\big\|C\big(x_1^{(m)}-x_1\big)\big\|_{E_1} &\leq a_1\big\|x_1^{(m)}-x_1\big\|_{E_1}+b_1\big\|A\big(x_1^{(m)}-x_1\big)\big\|_{E_1}\tolong 0,\\[2mm]
\big\|B\big(x_2^{(m)}-x_2\big)\big\|_{E_2} &\leq a_2\big\|x_2^{(m)}-x_2\big\|_{E_2}+b_2\big\|D\big(x_2^{(m)}-x_2\big)\big\|_{E_2}\tolong 0, 
\end{array}
\quad m\to\infty.
\end{align*}
Define $x^{(m)}:=\big(x_1^{(m)},x_2^{(m)}\big)\in\Phi, \,m\in\N$. The above convergences imply 
$x^{(m)}\to x$, $\cA x^{(m)}\to\cA x$ as $m\to\infty;$
hence $\Phi$ is a core of $\cA$.

Let $x:=(x_1,x_2)\in\Phi$. Then, for  $n\geq n_0(x):=\max\{n_1(x_1),n_2(x_2)\}$,
\begin{align*}
&\big(P_1^{(n)}\oplus P_2^{(n)}\big)x=\big(P_1^{(n)}x_1,P_2^{(n)}x_2\big)\in\dom(A^{(n)})\oplus\dom(D^{(n)})=\dom(\cA^{(n)}), \\
&\cA^{(n)}\big(P_1^{(n)}\oplus P_2^{(n)}\big)x- \cA x=\bmat \big(A^{(n)}P_1^{(n)}x_1-Ax_1\big) + \big(B^{(n)}P_2^{(n)}x_2-Bx_2\big)\\ \big(C^{(n)}P_1^{(n)}x_1-Cx_1\big)+ \big(D^{(n)}P_2^{(n)}x_2-Dx_2\big)\emat\tolong 0
\end{align*}
as $n\to\infty$. Now the claim follows from Theorem~\ref{propcoreforgsr}.
\end{proof}

\subsection{Results for bounded infinite operator matrices}\label{subsectionconvlemmas}
In this subsection we prove some useful lemmas about convergence of bounded infinite operator matrices;
the unbounded case is considered in the next subsection. 

\begin{definition}
Let $E_j^{(0)}$, $j\in\N$, be Banach spaces.  Define \index{$l^2\left(E_{j}: j\in\N\right)$}
\begin{align*}
\cE^{(0)}&:=l^2\big(E_j^{(0)}:j\in\N\big):=\bigg\{(x_j)_{j\in\N}:\, x_j\in E_j^{(0)},\,\sum_{j\in\N}\|x_j\|_{E_j^{(0)}}^2<\infty\bigg\},\\
\|x\|_{\cE^{(0)}}&:=\bigg(\sum_{j\in\N}\|x_j\|_{E_j^{(0)}}^2\bigg)^{\frac{1}{2}}, \quad x=(x_j)_{j\in\N}\in\cE^{(0)}.
\end{align*}
For each $j\in\N$ let $E_j,\,E_j^{(n)}\subset E_j^{(0)}, \,n\in\N,$ be closed subspaces;
let  $P_j:E_j^{(0)}\to E_{j}$ and $P_j^{(n)}:E_j^{(0)}\to E_j^{(n)}, \,n\in\N,$ be projections onto the respective subspaces.
Then  $$\cE:=l^2\left(E_{j}: j\in\N\right), \quad \cE^{(n)}:=l^2\big(E_j^{(n)}: j\in\N\big)\subset \cE^{(0)}, \quad n\in\N,$$ are Banach spaces. 
Denote the projections of $\cE^{(0)}$ onto the respective subspaces
by the diagonal block operator matrices
\begin{align*} \cP&:=\diag(P_{j}:\,j\in\N), \quad \cP^{(n)}:=\diag\big(P_j^{(n)}:\,j\in\N\big), \quad n\in\N.\end{align*}
\end{definition}

The following lemma is very useful for applications. It may be viewed as a Lebesgue's dominated convergence theorem  with respect to a counting measure. 

\begin{lemma} \label{domconv}
Let $x^{(n)}:=\big(x_j^{(n)}\big)_{j\in\N}\in \cE^{(0)}$, $n\in\N,$
and $x:=(x_j)_{j\in\N}\in\cE^{(0)}$ with
\beq \forall\,j\in\N:\quad 
\big\|x_j^{(n)}\big\|_{E_j^{(0)}}\leq \|x_j\|_{E_j^{(0)}}, \quad n\in\N,\quad \big\|x_j^{(n)}\big\|_{E_j^{(0)}}\tolong 0, \quad n\to\infty.\label{eqconvofentries}\eeq 
Then $\big\|x^{(n)}\big\|_{\cE^{(0)}}\to 0, \,n\to\infty.$ 
\end{lemma}

\begin{notation}\label{defcBcBn}
For $i,j\in\N$, let $B_{ij}\in L(E_j,E_i)$ and $B_{ij}^{(n)}\in L\big(E_j^{(n)},E_i^{(n)}\big)$, $n\in\N$.\index{block operator matrix!infinite}
Consider the infinite block operator matrices $\cB$, $\cB^{(n)}$, $n\in\N$, in $\cE$, $\cE^{(n)}$, $n\in\N$, respectively, defined by $$ \cB:=(B_{ij})_{i,j=1}^{\infty}, \quad \cB^{(n)}:=\big(B_{ij}^{(n)}\big)_{i,j=1}^{\infty}, \quad n\in\N.$$
\begin{enumerate}[label=\rm{\roman{*})}]
\item For each $i\in\N$ denote by $N_i\in\N\cup\{\infty\}$ the set of indices 
$j$ such that $B_{ij}\neq 0$ or $B_{ij}^{(n)}\neq 0$ for some $n\in\N$.
\item
For each $j\in\N$ denote by $M_j\in\N\cup\{\infty\}$ the set of indices 
$i$ such that $B_{ij}\neq 0$ or $B_{ij}^{(n)}\neq 0$ for some $n\in\N$.
\end{enumerate}
\end{notation}

For $\cB:=(B_{ij})_{i,j=1}^{\infty}\in L(\cE)$ and $\cB^{(n)}:=\big(B_{ij}^{(n)}\big)_{i,j=1}^{\infty}\in L(\cE^{(n)})$, $n\in\N$, we consider the following cases:
\begin{enumerate}[label=\rm{(\alph{*})}]
\item \label{casetridiag}
We have $N:=\sup_{i\in\N} \#N_i<\infty$, $M:=\sup_{j\in\N} \#M_j<\infty,$
and there exists $C\geq 0$ such that 
$$\forall\,i,j\in\N:\quad \big\| B_{ij}\big\|\leq C, \quad \big\|B_{ij}^{(n)}\big\|\leq C, \quad n\in\N.$$

\item \label{casenontridiag}
There exist $C_i\geq 0, \,i\in\N,$ such that 
$\sum_{i=1}^{\infty}C_i^2 \#N_i<\infty$
and
$$\forall\,i,j\in\N:\quad  \big\| B_{ij}\big\|\leq C_i, \quad \big\|B_{ij}^{(n)}\big\|\leq C_i, \quad n\in\N.$$

\item \label{casenontridiag2}
There exist $D_j\geq 0, \,j\in\N,$ such that 
$\sum_{j=1}^{\infty}D_j^2 \#M_j<\infty$ and
$$\forall\,i,j\in\N:\quad \big\| B_{ij}\big\|\leq D_j, \quad \big\|B_{ij}^{(n)}\big\|\leq D_j, \quad n\in\N.$$
\end{enumerate}

\begin{rem}
\begin{enumerate}[label=\rm{\roman{*})}]
\item Typical examples for case~\ref{casetridiag} are (uniformly) \emph{banded operator matrices} \index{block operator matrix!(semi-) banded}
for which there exists $L\in\N$ with $$B_{ij}=B_{ij}^{(n)}=0, \quad |i-j|>L, \quad n\in\N;$$ 
in this case we have $N\leq 2L+1$, $M\leq 2L+1$.

\item Typical examples for cases~\ref{casenontridiag} and \ref{casenontridiag2} are respectively lower and upper (uniformly) \emph{semibanded operator matrices}
for which there exists $L\in\N$ with 
$$B_{ij}=B_{ij}^{(n)}=0 \quad\text{ for } \begin{cases}j-i>L &\text{(case \ref{casenontridiag})},\\ i-j>L &\text{(case \ref{casenontridiag2})};\end{cases}$$
 we then have $\# N_i\leq i+L$ and $\# M_j\leq j+L$, respectively.

\item Note that, in general, neither \ref{casetridiag} is a special case of  \ref{casenontridiag} or \ref{casenontridiag2} nor vice versa.
In fact, for \ref{casenontridiag} and \ref{casenontridiag2},
 the norms of the entries need to decrease as $i\to\infty$ (in \ref{casenontridiag}) or $j\to\infty$ (in~\ref{casenontridiag2}), 
whereas in case \ref{casetridiag} the norms just need to be uniformly bounded;
vice versa, for \ref{casetridiag}, $(\# N_i)_{i\in\N}$ and $(\# M_j)_{j\in\N}$ need to be bounded sequences, whereas for \ref{casenontridiag} an \ref{casenontridiag2} they may be unbounded.
\end{enumerate}
\end{rem}

\begin{prop} \label{thmconvbddl2matrix}
Let $\cB:=(B_{ij})_{i,j=1}^{\infty}\in L(\cE)$ and $\cB^{(n)}:=\big(B_{ij}^{(n)}\big)_{i,j=1}^{\infty}\in L(\cE^{(n)})$, $n\in\N$, satisfy one of the cases {\rm (a), (b), (c)}.
If 
\beq \forall\,i,j\in\N:\quad B_{ij}^{(n)}P_j^{(n)}\slong B_{ij}P_j, \quad n\to\infty, \label{eqconvforiandl}\eeq
then we have 
$\cB^{(n)}\cP^{(n)}\slong\cB\cP$ as $n\to\infty.$
\end{prop}

\begin{proof}
(a) Let $y:=(y_j)_{j\in\N}\in \cE^{(0)}$.
For each $i\in\N$ choose $e_i\in E_i^{(0)}$ such that $\|e_i\|_{E_i^{(0)}}=1$ and define, with the constant $C$ from~\ref{casetridiag},
\begin{align*}
x_i&:=\bigg(\sum_{j\in N_i} \|y_j\|_{E_j^{(0)}}\bigg)2Ce_i,\quad
 x_i^{(n)}:=\sum_{j\in N_i}\left(B_{ij}^{(n)}P_j^{(n)}- B_{ij}P_j\right)y_j\in E_i^{(0)},\quad n\in\N.
\end{align*}
The element $x:=(x_i)_{i\in\N}$ belongs to $\cE^{(0)}$ since 
\begin{align*}
\|x\|_{\cE^{(0)}}^2&=
\sum_{i=1}^{\infty}\|x_i\|_{E_i^{(0)}}^2
= 4C^2\sum_{i=1}^{\infty}\bigg(\sum_{j\in N_i} \|y_j\|_{E_j^{(0)}}\bigg)^2 
\leq 4C^2\sum_{i=1}^{\infty}\#N_i \sum_{j\in N_i} \|y_j\|_{E_j^{(0)}}^2\\
&\leq  4C^2 N  \sum_{j=1}^{\infty}M_j\|y_j\|_{E_j^{(0)}}^2
\leq 4C^2N M\|y\|_{\cE^{(0)}}^2<\infty.
\end{align*}
We fix an $i\in\N$. Since $N_i$ is a finite set, the convergences in~\eqref{eqconvforiandl} imply 
$$\big\|x_i^{(n)}\big\|_{E_i^{(0)}}\leq \big\|x_i\big\|_{E_i^{(0)}}, \quad n\in\N,\quad \big\|x_i^{(n)}\big\|_{E_i^{(0)}}\tolong 0, \quad n\to\infty.$$
Now Lemma~\ref{domconv} applied to $x^{(n)}:=\big(x_i^{(n)}\big)_{i\in\N}\in\cE^{(0)}, \,n\in\N,$ yields
\be\big\|\big(\cB^{(n)}\cP^{(n)}-\cB\cP\big)y\big\|_{\cE^{(0)}}=\big\|x^{(n)}\big\|_{\cE^{(0)}}\tolong 0, \quad n\to\infty.\ee

(b) Proceed as in \ref{casetridiag} with $C$ replaced by $C_i$ in the definition of $x_i$. Then
\begin{align*}
\sum_{i=1}^{\infty}\|x_i\|_{E_i^{(0)}}^2
&\leq 4\sum_{i=1}^{\infty}C_i^2\Bigg(\sum_{j\in N_i} \|y_j\|_{E_j^{(0)}}\Bigg)^2
\leq  4\sum_{i=1}^{\infty}C_i^2 \# N_i \sum_{j=1}^{\infty}\|y_j\|_{E_j^{(0)}}^2<\infty.
\end{align*}
It is left to show 
\beq \forall\,i\in\N:\,\big\|x_i^{(n)}\big\|_{E_i^{(0)}}\tolong 0, \quad n\to\infty;\label{eqconvxin}\eeq
then the rest of the proof is completely analogous to \ref{casetridiag}.
We fix an $i\in\N$ and let $\eps>0$ be arbitrary. 
First note that 
\begin{align*}\big\|x_i^{(n)}\big\|_{E_i^{(0)}}&\leq \sum_{j=1}^{\infty}\underbrace{\big\|\big(B_{ij}^{(n)}P_j^{(n)}- B_{ij}P_j\big)y_j\big\|_{E_j^{(0)}}}_{\leq 2C_i\|y_j\|_{E_j^{(0)}}}, \quad n\in\N, 
\quad \sum_{j=1}^{\infty}2C_i\|y_j\|_{E_j^{(0)}} =\|x_i\|_{E_i^{(0)}}.
\end{align*}
There exists $j_{\eps}\in\N$ such that $$  \sum_{j=j_{\eps}}^{\infty}\big\|\big(B_{ij}^{(n)}P_j^{(n)}- B_{ij}P_j\big)y_j\big\|_{E_j^{(0)}}\leq \sum_{j=j_{\eps}}^{\infty} 2C_i \|y_j\|_{E_j^{(0)}}<\frac{\eps}{2}, \quad n\in\N.$$
The convergences in~\eqref{eqconvforiandl} imply the existence of $N_{\eps}\in\N$ such that $$ \sum_{j=1}^{j_{\eps}-1}\big\|\big(B_{ij}^{(n)}P_j^{(n)}- B_{ij}P_j\big)y_j\big\|_{E_j^{(0)}}<\frac{\eps}{2}, \quad n\geq N_{\eps}.$$
Altogether we have $\big\|x_i^{(n)}\big\|_{E_i^{(0)}}<\eps$ for all $n\geq N_{\eps}$; hence~\eqref{eqconvxin} is satisfied.

(c) For $i\in\N$ define $x_i^{(n)}$, $n\in\N$, as in \ref{casetridiag}, and set 
$$x_i:=\bigg(\sum_{j\in N_i} D_j \|y_j\|_{E_j^{(0)}}\bigg)2e_i\in E_i^{(0)}, \quad i\in\N.$$
Then $x:=(x_i)_{i\in\N}$ belongs to $\cE^{(0)}$ since 
\begin{align*}
\sum_{i=1}^{\infty}\|x_i\|_{E_i^{(0)}}^2& 
\leq 4 \sum_{i=1}^{\infty}\sum_{j\in N_i} D_j^2\sum_{k\in N_i} \|y_{k}\|_{E_{k}^{(0)}}^2
\leq 4 \sum_{j=1}^{\infty} D_j^2 (\# M_j) \sum_{k=1}^{\infty} \|y_{k}\|_{E_{k}^{(0)}}^2<\infty.
\end{align*}
The rest of the proof is analogous to \ref{casenontridiag}.
\end{proof}

\begin{corollary} \label{corcPnscP}
If $P_j^{(n)}\s P_j, \,n\to\infty,$ for all $j\in\N$, then $\cP^{(n)}\s\cP, \, n\to\infty$.
\end{corollary}

\begin{proof}
The claim follows immediately from Proposition~\ref{thmconvbddl2matrix}, case~\ref{casetridiag}.
\end{proof}

\subsection{Results for unbounded infinite operator matrices}\label{subsectioninfititegsr}

In this subsection we analyze whether a sequence of diagonally dominant infinite operator matrices converges in gsr-sense if the sequences for all diagonal elements do so.

Denote by  $\cE^{(0)}, \cE, \cE^{(n)}, \,n\in\N,$ the same spaces as in the previous subsection.
For $i,j\in\N$, let $A_{ij}:E_j\to E_i$ and $A_{ij}^{(n)}:E_j^{(n)}\to E_i^{(n)}$, $n\in\N$, be closable and densely defined,\index{block operator matrix!infinite}
and let
\begin{alignat*}{3}
 \cA&:=(A_{ij})_{i,j=1}^{\infty}, &\quad \dom(\cA)&:=l^2\Big(\underset{i\in\N}{\bigcap}\dom(A_{ij}):\,j\in\N\Big),\\
 \cA^{(n)}&:=\big(A_{ij}^{(n)}\big)_{i,j=1}^{\infty}, &\quad \dom\big(\cA^{(n)}\big)&:=l^2\Big(\underset{i\in\N}{\bigcap}\dom\big(A_{ij}^{(n)}\big):\,j\in\N\Big), \quad n\in\N.
\end{alignat*}

\noindent
We assume that $\cA$ and $\cA^{(n)}$, $n\in\N$, are densely defined in $\cE$ and $\cE^{(n)}$, $n\in\N$, respectively.

\begin{definition}\label{defdiagdominantinfinite}\index{block operator matrix!infinite diagonally dominant}
The block operator matrix $\cA$ is called
\emph{diagonally dominant}  if, for every $j\in\N$, the operators $A_{ij}$, $i\in\N$, are $A_{jj}$-bounded. 
Diagonal dominance of $\cA^{(n)}$, $n\in\N$, is defined analogously.
\end{definition}

\begin{theorem} \label{gsrdiag}
Assume that $\cA$, $\cA^{(n)}$, $n\in\N$, are diagonally dominant.
Let
\beq \label{defcTcTn}
\begin{aligned}
 \cT&:=\diag(A_{jj}:\,j\in\N),  \quad &\dom(\cT)&\!:=l^2\big(\dom(A_{jj}):\,j\in\N\big), \\
\cT^{(n)}&:=\diag\big(A_{jj}^{(n)}:\,j\in\N\big),   \quad &\dom\big(\cT^{(n)}\big)&\!:=l^2\big(\dom\big(A_{jj}^{(n)}\big):\,j\in\N\big),\quad n\in\N,
\end{aligned}
\eeq
and
\begin{align*}
\cS:=\cA-\cT, \quad \cS^{(n)}:=\cA^{(n)}-\cT^{(n)}, \quad n\in\N.
 \end{align*}
Suppose that there exist $\lm\in\underset{n\in\N}{\bigcap}\rho\big(\cT^{(n)}\big)\cap\rho(\cT)$ and $C_{\lm}>0$, $\gamma_{\lm}<1$ such that
\begin{align}
\begin{aligned}
\big\|(\cT^{(n)}-\lm)^{-1}\big\|&\leq C_{\lm}, \quad n\in\N,\\
\left\|\cS(\cT-\lm)^{-1}\right\|&<1,\quad \big\|\cS^{(n)}(\cT^{(n)}-\lm)^{-1}\big\|\leq \gamma_{\lm},\quad n\in\N.\label{eqSTbddinfinite}
\end{aligned}
\end{align}
Further assume that the bounded operator matrices $$\cB:=\cS(\cT-\lm)^{-1}, \quad \cB^{(n)}:=\cS^{(n)}(\cT^{(n)}-\lm)^{-1}, \quad n\in\N,$$ satisfy one of the cases~{\rm(a), (b), (c)} after Notation{\rm~\ref{defcBcBn}}.
If
 \begin{align*}
\begin{array}{lrll}
 \forall\,j\in\N:\quad &(A_{jj}^{(n)}-\lm)^{-1}P_j^{(n)}&\slong &(A_{jj}-\lm)^{-1}P_j,\\[1mm]
\forall\,i,j\in\N:\quad &A_{ij}^{(n)}(A_{jj}^{(n)}-\lm)^{-1}P_j^{(n)}&\slong &A_{ij}(A_{jj}-\lm)^{-1}P_j,\\[1mm]
\forall\,j\in\N:\quad &P_j^{(n)}&\slong &P_j,
 \end{array} \qquad n\to\infty,\label{eqconvforalldiag}
 \end{align*}
then $\lm\in\underset{n\in\N}{\bigcap}\rho(\cA^{(n)})\cap\rho(\cA)$ and 
$ (\cA^{(n)}-\lm)^{-1}\cP^{(n)}\s(\cA-\lm)^{-1}\cP$ as $n\to\infty.$
\end{theorem}

\begin{rem}
 The inequalities~\eqref{eqSTbddinfinite} imply that $\cS$ is $\cT$-bounded with $\cT$-bound~$<1$ and $\cS^{(n)}$ is $\cT^{(n)}$-bounded with $\cT^{(n)}$-bound $\leq\gamma_{\lm}<1$. 
\end{rem}

\begin{proof}[Proof of Theorem~{\rm\ref{gsrdiag}}]
Let $\lm\in\underset{n\in\N}{\bigcap}\rho\big(\cT^{(n)}\big)\cap\rho(\cT)$, $C_{\lm}>0$ and $\gamma_{\lm}<1$ satisfy the assumptions, and set $C:=\max\{\|(\cT-\lm)^{-1}\|,C_{\lm}\}$.
Then
\begin{align*}
\left\| (A_{jj}-\lm)^{-1}\right\|&\leq \left\| (\cT-\lm)^{-1}\right\|\leq C, \\
\big\|\big(A_{jj}^{(n)}-\lm\big)^{-1}\big\|&\leq \big\|(\cT^{(n)}-\lm)^{-1}\big\|\leq C, \quad n\in\N.
 \end{align*}
Proposition~\ref{thmconvbddl2matrix}, case~\ref{casetridiag} applied to $\left(\cT-\lm\right)^{-1}$ and $\left(\cT^{(n)}-\lm\right)^{-1}, \,n\in\N,$ yields
$$(\cT^{(n)}-\lm)^{-1}\cP^{(n)}\slong(\cT-\lm)^{-1}\cP, \quad n\to\infty.$$
By applying Proposition~\ref{thmconvbddl2matrix} to $\cB$, $\cB^{(n)}$, $n\in \N$, we obtain, for all of the cases (a),~(b),~(c), 
$$\cS^{(n)}(\cT^{(n)}-\lm)^{-1}\cP^{(n)}\slong \cS\left(\cT-\lm\right)^{-1}\cP, \quad n\to\infty.$$
Corollary~\ref{corcPnscP} yields $\cP^{(n)}\s\cP, \, n\to\infty$. 
Now the assertion follows from Theorem~\ref{perturbgsr}.
\end{proof}

The following example illustrates Theorem~\ref{gsrdiag}, case~(c).
\begin{example}
Let $\cE^{(0)}=\cE=\cE^{(n)}:=l^2(\N)$, $n\in\N$.
We define upper triangular matrices $\cA:=(A_{ij})_{i,j=1}^{\infty}, \cA^{(n)}:=\big(A_{ij}^{(n)}\big)_{i,j=1}^{\infty}\in C(l^2(\N))$ by
\begin{align*}
A_{ij}:=\begin{cases} j, &i< j,\\ j^3, &i=j,\\ 0, &\text{otherwise},\end{cases}
\quad 
A_{ij}^{(n)}:=\begin{cases} j, &i< j \leq n,\\ j^3, &i=j\leq n,\\ 0, &\text{otherwise}.\end{cases}
\end{align*}
Let $\cT$, $\cT^{(n)}$, $n\in\N$, be defined as in~\eqref{defcTcTn}.
Since the latter are selfadjoint operators,
it is obvious that $\C\backslash\R\subset\underset{n\in\N}{\bigcap}\rho\big(\cT^{(n)}\big)\cap\rho(\cT)$
and, for every $\lm\in\C\backslash\R$,
$$\forall\,j\in\N:\quad \big(A_{jj}^{(n)}-\lm\big)^{-1}\slong (A_{jj}-\lm)^{-1}, \quad n\to\infty.$$
Moreover, we have, for $i\neq j$,
\begin{align*}
A_{ij}(A_{jj}-\lm)^{-1}=\begin{cases}\frac{ j}{j^3-\lm}, &i< j,\\ 0, &\text{otherwise},\end{cases}
\quad 
A_{ij}^{(n)}\big(A_{jj}^{(n)}-\lm\big)^{-1}=\begin{cases} \frac{j}{j^3-\lm}, &i<j \leq n,\\ 0, &\text{otherwise}.\end{cases}
\end{align*}
Obviously, for $|\im(\lm)|$ sufficiently large, the matrices satisfy the assumptions  of  Theorem~\ref{gsrdiag}, case~(c),
with $$D_j:=\frac{j}{j^3-\lm}=\frac{1}{j^2-\frac{\lm}{j}}, \quad \# M_j=j-1, \quad j\in\N.$$
So we conclude $\lm\in\underset{n\in\N}{\bigcap}\rho(\cA^{(n)})\cap\rho(\cA)$ and 
$(\cA^{(n)}-\lm)^{-1}\s(\cA-\lm)^{-1}$ as $n\to\infty.$
\end{example}

\section{Discretely compact resolvents}\label{sectionsuffcondforcompres}
In this section we establish sufficient conditions for  a sequence of operators $T_n$, $n\in\N$, in varying Banach spaces to have \emph{discretely compact resolvents},
i.e.\ $$\exists\,\lm_0\in\underset{n\in\N}{\bigcap}\rho(T_n):\quad \big((T_n-\lm_0)^{-1}\big)_{n\in\N}\text{ discretely compact};$$ \index{compactness!discretely compact resolvents}
\vspace{-4mm}

\noindent
for the definition of discrete compactness see Definition~\ref{defdiscreteconvofop}.

As in Subsection~\ref{subsectiondirectgsr}, we first prove direct criteria and perturbation results (see Subsection~\ref{subsectiondirectdiscretecomp}).
Secondly, we establish results for operator matrices, both $2\times 2$ and infinite matrices (see Subsection~\ref{subsectiondisccompmatrices}).
Finally, we study the case that the operator domains are contained in (varying) Sobolev spaces and derive discretely compact Sobolev embeddings (see Subsection~\ref{subsectionsobolev}).

Consider a Banach space $E_0$  and closed subspaces $E,\,E_n\subset E_0, \,n\in\N$; as usual, we assume that the corresponding  projections converge strongly, $P_n\s P$.

\subsection{Direct criteria and perturbation results}\label{subsectiondirectdiscretecomp}

The following result relates (discrete) resolvent compactness to (discrete) compactness of embeddings.

\begin{prop}
Let $T\in C(E)$ and $T_n\in C(E_n)$, $n\in\N$.
 Define the Banach spaces 
$D:=\dom(T)$, $D_n:=\dom(T_n)$, $n\in\N$, 
equipped with the
graph norms $\|\cdot\|_D:=\|\cdot\|_T$, $\|\cdot\|_{D_n}:=\|\cdot\|_{T_n}$, $n\in\N$, respectively.
 Let $J:D\to E, \,J_n:D_n\to E_n$, $n\in\N$, be the natural embeddings. 
\begin{enumerate}[label=\rm{\roman{*})}]
\item  Let $\rho(T)\neq\emptyset$. The operator $J$ is compact if and only if $T$ has compact resolvent. The analogous result holds for $J_n$ and $T_n$.
\item Let $\Delta_b\big((T_n)_{n\in\N}\big)\neq\emptyset$.
The sequence $(J_n)_{n\in\N}$ is discretely compact if and only if for some {\rm(}and hence for all{\rm)} $\lm\in\Delta_b\left((T_n)_{n\in\N}\right)$ there exists $n_0\in\N$ 
such that the sequence $\left((T_n-\lm)^{-1}\right)_{n\geq n_0}$ is discretely compact.
\end{enumerate}
\label{equivcomp}
\end{prop}

\begin{proof}
We prove claim ii); claim i) is similar and well-known.

Let $\lm\in\Delta_b\big((T_n)_{n\in\N}\big)$ and $n_0\in\N$ such that $\lm\in\rho(T_n)$ for $n\geq n_0$. Consider an infinite subset $I\subset\N$ and $y_n\in E_n$, $n\in I$, and define  $x_n:=(T-\lm)^{-1}y_n\in D_n$, $n\in I$.
 First note that  we have
\begin{align*}
\|x_n\|_{T_n}&=\|x_n\|_{E_n}+\|T_nx_n\|_{E_n}\leq (1+|\lm|)\|x_n\|_{E_n}+\|(T_n-\lm)x_n\|_{E_n} \\ &\leq \big((1+|\lm|)\|(T_n-\lm)^{-1}\|+1\big)\|y_n\|_{E_n},\\[1mm]
\|y_n\|_{E_n}&=\|(T_n-\lm)x_n\|_{E_n}\leq \|T_nx_i\|_{E_n}+|\lm|\,\|x_n\|_{E_n}\leq (1+|\lm|)\,\|x_n\|_{T_n}.
\end{align*}
Hence $(\|x_n\|_{T_n})_{n\in\N}$ is bounded if and only if $(\|y_n\|_{E_n})_{n\in I}$ is bounded.

Assume that $(J_n)_{n\in\N}$ is discretely compact and $(\|x_n\|_{T_n})_{n\in I}$ is bounded. 
Then the sequence of elements $(T_n-\lm)^{-1}y_n=x_n\in~D_n$, $n\in I,$ has a subsequence that converges in $E_0$ with limit in $E$; hence  $\left((T_n-\lm)^{-1}\right)_{n\geq n_0}$ is discretely compact.
Vice versa,  the discrete compactness of $\left((T_n-\lm)^{-1}\right)_{n\geq n_0}$ implies discrete compactness of $(J_n)_{n\in\N}$. 
\end{proof}

We prove a perturbation result for discrete compactness of the resolvents.
Note that the assumptions are similar to the ones used in the perturbation result for gsr-convergence (compare Theorem~\ref{perturbgsr}).

\begin{theorem}
Let $T_n\in C(E_n), \,n\in\N$.
Let  $S_n$, $n\in\N$, be linear operators in $E_n$, $n\in\N$, with $\dom(T_n)\subset\dom(S_n)$, $n\in\N$, respectively.
Define $$A_n:=T_n+S_n, \quad n\in\N.$$
Suppose that there exist $\lm\in\underset{n\in\N}{\bigcap}\rho(T_n)$ 
and $\gamma_{\lm}<1$ with
\begin{equation}\label{eqSTbdd2}
\left\|S_n(T_n-\lm)^{-1}\right\|\leq \gamma_{\lm},\quad n\in\N.
\end{equation}
Then $\lm\in\underset{n\in\N}{\bigcap}\rho(A_n)$.
If the sequence $\left((T_n-\lm)^{-1}\right)_{n\in\N}$ is discretely compact, then so is $\left((A_n-\lm)^{-1}\right)_{n\in\N}$.
\label{perturbdiscrcompres}
\end{theorem}

\begin{rem}
 The inequalities~\eqref{eqSTbdd2} imply that, for every $n\in\N$,  $S_n$ is $T_n$-bounded with $T_n$-bound $\leq\gamma_{\lm}<1$. 
\end{rem}

\begin{proof}[Proof of Theorem~{\rm\ref{perturbdiscrcompres}}]
 Let $\lm$ satisfy the assumptions. Then, using~\eqref{neumanngamma} and Lemma \ref{lemmapertdisccomp}~i),
the discrete compactness of  $\left((T_n-\lm)^{-1}\right)_{n\in\N}$ implies that the sequence  $\left((A_n-\lm)^{-1}\right)_{n\in\N}$ is discretely compact.
\end{proof}

The following result is an immediate consequence of Theorem \ref{perturbdiscrcompres} for the case that the perturbations are bounded operators.

\begin{corollary}
Let $T_n\in C(E_n), \,n\in\N$, and $S_n\in L(E_n)$, $n\in\N$.
Define $$A_n:=T_n+S_n, \quad n\in\N.$$
Suppose that there exist $\lm\in\underset{n\in\N}{\bigcap}\rho(T_n)$ 
and $\gamma_{\lm}<1$ with
\be \left\|S_n\right\|\leq \gamma_{\lm}\,\|(T_n-\lm)^{-1}\|^{-1},\quad n\in\N.\label{eqSSnrelbdd2}\ee
Then $\lm\in\underset{n\in\N}{\bigcap}\rho(A_n)$.
If the sequence $\left((T_n-\lm)^{-1}\right)_{n\in\N}$ is discretely compact, then so is $\left((A_n-\lm)^{-1}\right)_{n\in\N}$.
\label{corperturbdiscrcompres}
\end{corollary}


\subsection{Results for block operator matrices}\label{subsectiondisccompmatrices}
In this subsection we consider (finite and infinite)  operator matrices.
We study whether a sequence of diagonally dominant operator matrices 
$$\cA^{(n)}:=\big(T_{ij}^{(n)}\big)_{i,j=1}^{N}, \quad n\in\N,$$ 
has discretely compact resolvents if the sequences $\big(T_{ii}^{(n)}\big)_{n\in\N}$, $i\in\N$, of all diagonal entries have discretely compact resolvents.

First, we consider diagonally dominant $2\times 2$ operator matrices, i.e.\ $N=2$.
We use the same notation as in Subsection~\ref{subsection2x2matricesconv}.

\begin{theorem}
Suppose that the block operator matrices satisfy the following:
\begin{enumerate}[label=\rm{\roman{*})}]
\item the operators $A^{(n)}, \,n\in\N,$  and  $D^{(n)}, \,n\in\N,$ have discretely compact resolvents, respectively;
\item the matrices $\cA^{(n)}$, $n\in\N$, are diagonally dominant; 
\item 
there exist $\lm\in\underset{n\in\N}{\bigcap}\left(\rho(A^{(n)})\cap\rho(D^{(n)})\right)$ and constants $\gamma_{\lm}^{AC}, \,\gamma_{\lm}^{DB}\geq 0$ such that
$$\big((A^{(n)}-\lm)^{-1}\big)_{n\in\N}, \quad \big((D^{(n)}-\lm)^{-1}\big)_{n\in\N}$$ are bounded sequences,
\beq \big\|C^{(n)}(A^{(n)}-\lm)^{-1}\big\|\leq\gamma_{\lm}^{AC},\quad \big\|B^{(n)}(D^{(n)}-\lm)^{-1}\big\|\leq\gamma_{\lm}^{DB}, \quad n\in\N,\label{ineqgammalmforACandDB}\eeq
and \beq \gamma_{\lm}^{AC} \gamma_{\lm}^{DB}<1.\label{eqgammas}\eeq
\end{enumerate}
Then 
$\lm\in\underset{n\in\N}{\bigcap}\rho(\cA^{(n)})$ and the sequence
 $\left((\cA^{(n)}-\lm)^{-1}\right)_{n\in\N}$ is discretely compact.
\label{propmatrixwithdiscrcompres}
\end{theorem}

\begin{rem}\label{remgammaorder}
 The inequalities~\eqref{ineqgammalmforACandDB} imply that $\cA^{(n)}$ is diagonally dominant of order $\max\{\gamma_{\lm}^{AC},\gamma_{\lm}^{DB}\}$ (see Definition~\ref{defdiagdominant}~i)).
Note that~\eqref{eqgammas} is more general than $\max\{\gamma_{\lm}^{AC},\gamma_{\lm}^{DB}\}<1$
and it implies $\min\{\gamma_{\lm}^{AC},\gamma_{\lm}^{DB}\}<1$.
\end{rem}

\begin{proof}[Proof of Theorem~{\rm\ref{propmatrixwithdiscrcompres}}]
Let $\lm$  satisfy the assumption~iii).
We show the claim in the case $\gamma_{\lm}^{DB}=\min\big\{\gamma_{\lm}^{AC},\gamma_{\lm}^{DB}\big\}<1;$ the other case is analogous.
For each $n\in\N$ define $$\cT^{(n)}:=\bmat A^{(n)} & 0\\ C^{(n)} &D^{(n)}\emat, \quad \cS^{(n)}:=\cA^{(n)}-\cT^{(n)}=\bmat 0 & B^{(n)}\\ 0 & 0\emat.$$
Then, for each $n\in\N$, 
\begin{align*}
(\cT^{(n)}-\lm)^{-1}=&\bmat (A^{(n)}-\lm)^{-1} & 0\\ 0 & 0\emat+ \bmat 0 & 0\\ 0 & (D^{(n)}-\lm)^{-1}\emat\\[1mm]
&+\bmat 0 & 0 \\ -(D^{(n)}-\lm)^{-1}C^{(n)}(A^{(n)}-\lm)^{-1} & 0\emat;
\end{align*}
denote the matrices on the right-hand side by $\mathcal U^{(n)}$, $\mathcal V^{(n)}$, $\mathcal W^{(n)}$, respectively.
The assumption~i) implies that $\big(\mathcal U^{(n)}\big)_{n\in\N}$, $\big(\mathcal V^{(n)}\big)_{n\in\N}$ are discretely compact sequences.
By the assumption~iii), $\big(C^{(n)}(A^{(n)}-\lm)^{-1}\big)_{n\in\N}$ is a bounded sequence. Now i) and Lemma~\ref{lemmapertdisccomp}~i) imply the discrete compactness of $\big(\mathcal W^{(n)}\big)_{n\in\N}$.
By Lemma~\ref{lemmapertdisccomp}~iii), the sequence $\big((\cT^{(n)}-\lm)^{-1}\big)_{n\in\N}$ is discretely compact.
Assumption iii) yields the estimate
\begin{align*}
&\big\|\cS^{(n)}(\cT^{(n)}-\lm)^{-1}\big\|\\
 &\leq\max\left\{\big\|B^{(n)}(D^{(n)}-\lm)^{-1}C^{(n)}(A^{(n)}-\lm)^{-1}\big\|, \,\big\|B^{(n)}(D^{(n)}-\lm)^{-1}\big\|\right\}\\
&\leq \max\{\gamma_{\lm}^{AC} \gamma_{\lm}^{DB},\,\gamma_{\lm}^{DB}\}<1, \quad n\in\N.
\end{align*}
Now the claim follows from Theorem~\ref{perturbdiscrcompres}.
\end{proof}

\begin{example}
Let the operators $A^{(n)}, \,n\in\N,$  and  $D^{(n)}, \,n\in\N,$ be selfadjoint and have discretely compact resolvents, respectively.
If the operators $B^{(n)}, C^{(n)},\,n\in\N,$ are uniformly bounded, then the sequence of operators $$\cA^{(n)}:=\bmat A^{(n)} &B^{(n)}\\ C^{(n)} &D^{(n)}\emat, \quad n\in\N,$$ has discretely compact resolvents.
This follows from Theorem~\ref{propmatrixwithdiscrcompres} since for every $\gamma<1$ there exists 
$\lm\in\C\backslash\R$  with $|\im(\lm)|$ sufficiently large such that the assumption~iii) of Theorem~\ref{propmatrixwithdiscrcompres} holds with 
$$\gamma_{\lm}^{AC}=\gamma_{\lm}^{DB}:=\frac{1}{|\im(\lm)|}\sup_{n\in\N}\,\max\big\{\|B^{(n)}\|,\|C^{(n)}\|\big\}.$$
\end{example}

Now we  study sequences of diagonal operator matrices; first for finitely many diagonal elements, then for infinitely many.
For block operator matrices that are not diagonal, one may combine Theorem~\ref{discrcompresdiagfinite}/\ref{discrcompresdiag} 
with Theorem~\ref{perturbdiscrcompres} or its corollary (as in Theorem~\ref{propmatrixwithdiscrcompres} for the $2\times 2$ case).

Let $E_j^{(0)}$, $j\in\N$, be Banach spaces and, for every $j\in\N$, let  $E_j,\,E_j^{(n)}\subset E_j^{(0)}$, $n\in\N$, be closed subspaces.

\begin{theorem}
 Let $k\in\N$ and consider Banach spaces $$\cE:=\underset{j=1,\dots,k}{\bigoplus}E_j, \quad \cE^{(n)}:=\underset{j=1,\dots,k}{\bigoplus}E_j^{(n)}\subset \cE^{(0)}:=\underset{j=1,\dots,k}{\bigoplus}E_j^{(0)}, \quad n\in\N.$$ 
Define
$$\cT^{(n)}:=\diag\big(T_j^{(n)}:\,j=1,\dots, k\big)\in C(\cE^{(n)}), \quad n\in\N. $$
Suppose that there exists $\lm\in\underset{j=1,\dots,k}{\bigcap}\,\,\underset{n\in\N}{\bigcap}\rho\big(T_j^{(n)}\big)$ such that 
$$\big((T_j^{(n)}-\lm)^{-1}\big)_{n\in\N}, \quad j=1,\dots,k,$$ are discretely compact sequences.
Then we have $\lm\in\underset{n\in\N}{\bigcap}\rho(\cT^{(n)})$ and the sequence $\big((\cT^{(n)}-\lm)^{-1}\big)_{n\in\N}$
is discretely compact.
\label{discrcompresdiagfinite}
\end{theorem}

\begin{proof}
For $k=2$ the claim is an immediate consequence of Lemma~\ref{lemmapertdisccomp}~iii) applied to 
\be A_n^{(1)}=\bmat \big(T_1^{(n)}-\lm\big)^{-1} & 0 \\ 0 & 0\emat, \quad A_n^{(2)}=\bmat 0 & 0 \\ 0 & \big(T_2^{(n)}-\lm\big)^{-1}\emat, \quad n\in\N.\ee
For $k\in\N$ with $k>2$ the claim follows by induction.
\end{proof}

\begin{theorem}
Assume that $E_j^{(0)}=E_j$, $j\in\N$.
Consider  Banach spaces $$\cE^{(n)}:=l^2\big(E_j^{(n)}: j\in\N\big)\subset \cE:=l^2\big(E_j:j\in\N\big), \quad n\in\N.$$  
Define $$\cT^{(n)}:=\diag\big(T_j^{(n)}:\,j\in\N\big)\in C(\cE^{(n)}), \quad n\in\N.$$
 Suppose that there exists $\lm\in\underset{j\in\N}{\bigcap}\,\,\underset{n\in\N}{\bigcap}\rho\big(T_j^{(n)}\big)$ such that 
$$\big((T_j^{(n)}-\lm)^{-1}\big)_{n\in\N}, \quad j\in\N,$$ are discretely compact sequences.
We further assume that 
\beq \label{univresconv}
\sup_{n\in\N}\big\|\big(T_j^{(n)}-\lm\big)^{-1}\big\|\tolong 0, \quad j\to\infty.
\eeq
Then we have $\lm\in\underset{n\in\N}{\bigcap}\rho(\cT^{(n)})$ and the sequence  $\big((\cT^{(n)}-\lm)^{-1}\big)_{n\in\N}$
is discretely compact.
\label{discrcompresdiag}
\end{theorem}

\begin{proof}
Let $\lm$ satisfy the assumptions. Fix an $n\in\N$. Since, by the  assumption~\eqref{univresconv}, the sequence  $\big((T_j^{(n)}-\lm)^{-1}\big)_{j\in\N}$ is bounded,  $\lm$ belongs to $\rho(\cT^{(n)})$.
We define
\begin{align*}
\cA^{(n)}&:=(\cT^{(n)}-\lm)^{-1}=\diag\big((T_j^{(n)}-\lm)^{-1}:\,j\in\N\big), \\
\cA^{(n;k)}&:=\diag\big((T_j^{(n)}-\lm)^{-1}:\,j=1,\dots,k\big)\oplus 0\in L\big(\cE^{(n)}\big), \quad k\in\N.
\end{align*}
From the assumption~\eqref{univresconv}, it follows that
$$\lim_{k\to\infty}\,\sup_{n\in\N}\big\|\cA^{(n)}-\cA^{(n;k)}\big\|=\lim_{k\to\infty}\,\sup_{j>k}\,\sup_{n\in\N}\big\|(T_j^{(n)}-\lm)^{-1}\big\|=0.$$
Since, by the assumptions, $\big((T_j^{(n)}-\lm)^{-1}\big)_{n\in\N}, \, j\in\N,$ are discretely compact sequences, Theorem  \ref{discrcompresdiagfinite} implies that
$\big(\cA^{(n;k)}\big)_{n\in\N}$ is discretely compact for each $k\in\N$. Altogether, Proposition \ref{propapproxdiscrcomp} yields that $(\cA^{(n)})_{n\in\N}=\big((\cT^{(n)}-\lm)^{-1}\big)_{n\in\N}$ is discretely compact.
\end{proof}


\subsection{Discretely compact Sobolev embeddings}\label{subsectionsobolev}

In this subsection we assume that $E_0:=L^p(\R^d)$ and $E:=L^p(\Omega)$, where  $p\in(1,\infty)$, $d\in\N$ and $\Omega\subset\R^d$ is an open subset.
The space $E_n$ is assumed to be $L^p(\Omega_n)$ for some open subset $\Omega_n\subset\R^d$ that may vary in $n\in\N$.
 Denote by ${\rm m}_d$\index{${\rm m}_d$}  the Lebesgue measure on $\R^d$.

We consider operators $T_n\in C(E_n), \,n\in\N,$ whose domains are assumed to be subspaces of Sobolev spaces $W^{m,p}(\Omega_n), \,n\in\N$, for some $m\in\N$.
In this case it is sufficient to study discrete compactness of the Sobolev embeddings $$\widetilde J_n: W^{m,p}(\Omega_n)\to L^p(\Omega_n), \quad n\in\N,$$ 
in order to conclude that the sequence of embeddings $J_n:D_n\to E_n$, $n\in\N$, (as defined in Proposition~\ref{equivcomp}) is discretely compact (see Theorem~\ref{thmsobolev} below).

\begin{definition}
 Let $\Omega,\, \Omega_n, \,n\in\N,$ be bounded open subsets of $\R^d$.
\begin{enumerate}[label=\rm{\roman{*})}]
\item The set $\Omega$ is said to have the \emph{segment property}\index{segment property} if  there exist a finite open covering   $\left\{O_k:\,k=1,\dots,r\right\}$ of $\partial\Omega$
and corresponding vectors 
 $\alpha_k\in\R^d\backslash\{0\}$, $k=1,\dots,r,$ such that 
$$\lbar{\Omega\cap O_k}+t\alpha_k\subset\Omega, \quad t\in(0,1),  \quad k=1,\dots,r.$$
\item The pair  $\left\{\Omega, \{\Omega_n:\,n\in\N\}\right\}$ is said to have the \emph{uniform segment property} if there exist an open covering   $\left\{O_k:\,k=1,\dots,r\right\}$ of $\partial\Omega$
and corresponding vectors 
$\alpha_k\in\R^d\backslash\{0\}, \,k=1,\dots,r,$ such that  $\left\{O_k:\,k=1,\dots,r\right\}$  is an open covering of $\partial\Omega_n$ for sufficiently large $n\in\N$, say $n\geq n_0$,
and
\begin{alignat*}{3} 
 \forall\,n\geq n_0:\quad &\lbar{\Omega_n\cap O_k}+t\alpha_k\subset\Omega_n, \quad t\in(0,1),  \quad&&k=1,\dots,r,\\
\forall\,\eps\in(0,1)\,\exists\,n_{\eps}\in\N:\quad &\lbar{\underset{n\geq n_{\eps}}{\bigcup}\left(\Omega_n\cap O_k\right)}+\eps\alpha_k\subset\Omega, \quad &&k=1,\dots,r.
\end{alignat*}
\end{enumerate}
\end{definition}

\begin{rem}\label{remsegmentprop}
\begin{enumerate}[label=\rm{\roman{*})}]
 \item 
It is easy to see that if the pair  $\left\{\Omega, \{\Omega_n:\,n\in\N\}\right\}$  has the uniform segment property, then $\Omega_n, \,n\geq n_0,$ all have the segment property.
\item 
If each compact subset $S\subset\Omega$ satisfies $S\subset\Omega_n$ for all sufficiently large $n\in\N$, then $\Omega\subset\bigcup_{n\geq n_{\eps}}\Omega_n$ for each $n_{\eps}\in\N$.
If, in addition, the pair  $\left\{\Omega, \{\Omega_n:\,n\in\N\}\right\}$  has the uniform segment property, then $\Omega$ has the segment property
since, for each $\eps\in (0,1)$, $$\lbar{\Omega\cap O_k}+\eps\alpha_k\subset\lbar{\underset{n\geq n_{\eps}}{\bigcup}\left(\Omega_n\cap O_k\right)}+\eps\alpha_k\subset\Omega,\quad k=1,\dots,r.$$
\end{enumerate}
\end{rem}

\begin{example}\label{examplesegment}
 \begin{enumerate}[label=\rm{\roman{*})}]
  \item The motivation for defining the segment property is that the interior of the set $\Omega$ shall not lie on both sides of the boundary. For instance, in $\R^2$ the set $B_1(0)\backslash \left((0,1)\times \{0\}\right)$ does not have the segment property.
  \item In dimension $d=1$, it is easy to see that an open bounded subset $\Omega\subset\R$ has the segment property if $\Omega$ is the finite union of open bounded intervals $I_l, \,l=1,\dots,L$, where two different intervals have positive distance.
Then $r=2L$ and each $O_k$ contains exactly one endpoint of an interval $I_l$; the number $\alpha_k$ is positive (negative) if it is the left (right) endpoint, with $|\alpha_k|$ less than the length of $I_l$.
 \item An example in dimension $d=1$ for bounded open sets $\Omega,\Omega_n\subset\R$, $n\in\N$, such that the pair $\left\{\Omega, \{\Omega_n:\,n\in\N\}\right\}$ has the uniform segment property is
\begin{alignat*}{3}
 \Omega&:=\underset{l=1,\dots,L}{\bigcup}(a_l,b_l), \quad &&a_l<b_l<a_{l+1},\\
 \Omega_n&:=\underset{l=1,\dots,L}{\bigcup}\big(a_l^{(n)},b_l^{(n)}\big), \quad &&a_l^{(n)}<b_l^{(n)}<a_{l+1}^{(n)},
\end{alignat*}
with $a_l^{(n)}\to a_l$, $b_l^{(n)}\to b_l$, $n\to\infty$, for every $l=1,\dots,L$.
 \end{enumerate}
\end{example}

In the following we establish discretely compact Sobolev embeddings using earlier results 
by Grigorieff~\cite{grigorieffsobolev}. 
As Grigorieff's results are confined to dimension $d\geq 2$, we prove the case $d=1$ separately; to avoid unnecessarily technicalities, in $d=1$ we prove the result directly for the case studied in 
Example~\ref{examplesegment}~iii) where $\Omega_n$, $n\in\N$, are unions of $L<\infty$ intervals whose endpoints converge to the ones of $\Omega$.

\begin{theorem}\label{thmsobolev}
For $d\geq 2$, 
suppose that $\Omega,\,\Omega_n\subset\R^d$, $n\in\N,$ are bounded open subsets that satisfy the following:
\begin{enumerate}[label=\rm{(\roman{*})}]
\item  each compact subset $S\subset\Omega$ is also a subset of $\Omega_n$ for all sufficiently large $n\in\N$;
\item the pair  $\left\{\Omega, \{\Omega_n:\,n\in\N\}\right\}$  has the uniform segment property;
\item we have ${\rm m}_d(\Omega_n\backslash\Omega)\to 0, \, n\to\infty.$ 
\end{enumerate}
For $d=1$,  suppose that $\Omega,\,\Omega_n\subset\R$, $n\in\N,$ are as in Example~{\rm\ref{examplesegment}~iii)}. 

 Let $T_n\in C(E_n), \,n\in\N$, with $\dom(T_n)\subset W^{m,p}(\Omega_n)$, $n\in\N$, for some $m\in\N$. 
If the embeddings $$B_n: \big(\dom(T_n),\|\cdot\|_{T_n}\big)\to W^{m,p}(\Omega_n), \quad n\in\N,$$
 are uniformly bounded, then the sequence $(J_n)_{n\in\N}$ of embeddings $$J_n: \big(\dom(T_n),\|\cdot\|_{T_n}\big)\to L^p(\Omega), \quad  n\in\N,$$ is discretely compact.
\end{theorem}

\begin{proof}
For dimension $d\geq 2$, the sequence $(\widetilde J_n)_{n\in\N}$ of embeddings
$$\widetilde J_n: W^{m,p}(\Omega_n)\to L^p(\Omega_n), \quad n\in\N,$$
 is discretely compact by \cite[Satz 4.(9)]{grigorieffsobolev}. 
Since $(B_n)_{n\in\N}$ is a bounded sequence by the assumptions,
the claim follows from Lemma~\ref{lemmapertdisccomp}~i). 

For dimension $d=1$, let $I\subset\N$ be an infinite subset and let $f_n\in\dom(T_n), \,n\in I,$ satisfy that $\left(\|f_n\|_{T_n}\right)_{n\in I}$ is bounded. Then $\left(\|f_n\|_{W^{m,p}(\Omega_n)}\right)_{n\in I}$ is bounded since $(B_n)_{n\in\N}$ is a bounded sequence.
 Define $$\Lambda:=\Omega\times(0,1), \quad \Lambda_n:=\Omega_n\times(0,1), \quad n\in\N.$$
These sets are bounded open  subsets of $\R^2$.
The idea of the proof is to show that $\Lambda, \Lambda_n$, $n\in\N,$ satisfy assumptions~(i)--(iii) for $d=2$;
then the sequence of embeddings 
\beq \widetilde J_n^{(2)}: W^{m,p}(\Lambda_n)\to L^p(\Lambda_n), \quad n\in\N,\label{eq.sobolevd2}\eeq
 is discretely compact by  \cite[Satz 4.(9)]{grigorieffsobolev}.
From this, at the end, we conclude that the sequence of elements $f_n\in L^p(\Omega_n)$, $n\in I$, has a convergent subsequence in $L^p(\R)$ with limit function in $L^p(\Omega)$.

It is easy to see that properties (i) and (iii) are satisfied for $\Lambda$, $\Lambda_n$, $n\in\N$.
It remains to check (ii), i.e.\ the uniform segment property. 
There exists $\delta>0$ such that $b_l-a_l>3\delta$ for all $l=1,\dots,L$ and $a_{l+1}-b_l>2\delta$ for all $l=1,\dots,L-1$.
Note that the latter implies $(b_l-\delta,b_l+\delta)\cap (a_{l+1}-\delta,a_{l+1}+\delta)=\emptyset$.
Since $a_l^{(n)}\to a_l$ and $b_l^{(n)}\to b_l$ as $n\to\infty$, there exists $n_0\in\N$
with $$\forall\,n\geq n_0:\quad \big|a_l^{(n)}-a_l\big|<\delta, \quad \big|b_l^{(n)}-b_l\big|<\delta, \quad l=1,\dots,L.$$
Since $b_l-a_l>3\delta$ we thus have, for all $t\in (0,1)$, 
\beq \forall\,n\geq n_0:\quad\big[a_l^{(n)},a_l+\delta\big]+t \delta\subset\Omega_n, \quad \big[b_l-\delta,b_l^{(n)}\big]-t \delta\subset\Omega_n, \quad l=1,\dots,L.\label{eqsegment1}\eeq
For any $\eps\in (0,1)$ there exists $n_{\eps}\geq n_0$ such that 
$$\sup_{n\geq n_{\eps}}\big|a_l^{(n)}-a_l\big|<\eps\delta, \quad \sup_{n\geq n_{\eps}}\big|b_l^{(n)}-b_l\big|<\eps\delta, \quad l=1,\dots,L.$$
Hence, again with $b_l-a_l>3\delta$,
\beq\label{eqsegment2}
\begin{array}{l}
 \overline{\bigcup_{n\geq n_{\eps}}\big(a_l^{(n)},a_l+\delta\big)}+\eps\delta\subset (a_l-\eps \delta,a_l+\delta]+\eps\delta\subset\Omega, \\[2mm]
 \overline{\bigcup_{n\geq n_{\eps}}\big(b_l-\delta,b_l^{(n)}\big)}+\eps\delta\subset [b_l-\delta,b_l+\eps\delta)-\eps\delta\subset\Omega, 
\end{array}\quad l=1,\dots,L.
\eeq
Now \eqref{eqsegment1} and \eqref{eqsegment2} imply that $\{\Lambda,\{\Lambda_n:\,n\in\N\}\}$
has the uniform segment property with finite open covering
\begin{align*}
\Big\{&(a_l-\delta,a_l+\delta)\times (-1/2,2/3),\,
(a_l-\delta,a_l+\delta)\times (1/3,3/2),\\
&(b_l-\delta,b_l+\delta)\times (-1/2,2/3),\,
(b_l-\delta,b_l+\delta)\times (1/3,3/2),\\
&(a_l+\delta/2,b_l-\delta/2)\times (-1/2,1/3),\,
(a_l+\delta/2,b_l-\delta/2)\times (2/3, 3/2):\,l=1,\dots,L\Big\}
\end{align*}
and corresponding set of non-zero vectors in $\R^2$
\begin{align*}
\bigg\{
\bmat \delta \\ 1/3\emat, \,
\bmat \delta\\ \!-1/3\emat,\,
\bmat -\delta\\ 1/3\emat,\,
\bmat -\delta\\ \!-1/3\emat,\,
\bmat 0\\ 1/3\emat,\,
\bmat 0\\ \!-1/3\emat: \,l=1,\dots,L\bigg\}.
\end{align*}

Altogether, assumptions~(i)--(iii) are satisfied and thus  \cite[Satz 4.(9)]{grigorieffsobolev} 
yields that  the sequence of embeddings 
in~\eqref{eq.sobolevd2} is discretely compact.
Note that $f_n\in W^{m,p}(\Lambda_n)$, $n\in I,$ and $\left(\|f_n\|_{W^{m,p}(\Lambda_n)}\right)_{n\in I}$ is bounded.
Hence there exist $f\in L^p(\Lambda)$ and an infinite subset $I_2\subset I$ such that $(\widetilde J_n^{(2)}f_{n})_{n\in I_2}$ converges to $f$ in $L^p(\R^2)$.
Since $f\in L^p(\Lambda)$,  we have $f(\cdot,x_2)\in L^p(\Omega)$ for almost all $x_2\in(0,1)$; denote by $\Theta_1\subset(0,1)$ the set of such $x_2$.
The convergence $\|f_n-f\|_{L^p(\R^2)}\to 0$ as $n\in I_2$, $n\to\infty$, 
implies the existence of an infinite subset $I_3\subset I_2$ so that, for $n\in I_3$,
$$\int_{\R}|f_{n}(x_1)-f(x_1,x_2)|^p\,\rd x_1\tolong 0, \quad n\to\infty,$$
 for almost all $x_2\in(0,1)$ (see e.g.\ \cite[Theorem~B.98~(iii)]{leoni} with $u\equiv 0, \,p=1$); denote by $\Theta_2\subset(0,1)$ the set of such $x_2$.
For $x_2\in\Theta_1\cap\Theta_2$ we hence obtain $f(\cdot,x_2)\in L^p(\Omega)$  and  $\|f_{n}-f(\cdot,x_2)\|_{L^p(\R)}\to 0$ as $n\in I_3$, $n\to\infty$.
So we have shown that $(J_nf_n)_{n\in I}$ has a convergent subsequence in $L^p(\R)$ with limit in $L^p(\Omega)$.
\end{proof}

\section{Applications to domain truncation method and Galerkin method}\label{sectionapplications}
In this section we give examples of spectrally exact operator approximations. All underlying spaces are Hilbert spaces $H$ and $H_n\subset H$, $n\in\N$, with corresponding orthogonal projections $P_n$ with $P_n\s I$.
For an operator $T\in C(H)$ and approximating operators $T_n\in C(H_n)$, $n\in\N$,  we check whether there exists an element
$\lm_0\in\underset{n\in\N}{\bigcap}\rho(T_n)\cap\rho(T)$ such that
\begin{enumerate}
\item[\rm (a)] $(T-\lm_0)^{-1}$, $(T_n-\lm_0)^{-1},\,n\in\N,$ are compact operators;
\item[\rm (b)] the sequence $\left((T_n-\lm_0)^{-1}\right)_{n\in\N}$ is discretely compact;
\item[\rm (c)] the resolvents converge strongly, $(T_n-\lm_0)^{-1}P_n\s (T-\lm_0)^{-1}$;
\item[\rm (d)] the adjoint resolvents converge strongly, $(T_n^*-\overline{\lm_0})^{-1}P_n\s (T^*-\overline{\lm_0})^{-1}$.
\end{enumerate}
If (a)--(c) are satisfied, then  Theorem~{\rm\ref{mainthmforgsr}} is applicable which yields, in particular, that $(T_n)_{n\in\N}$ is a spectrally exact approximation of $T$.
If, in addition, (d) holds, then Theorem~\ref{thmgsrimpliesgnrbasic} yields generalized norm resolvent convergence, i.e.\ $(T_n-\lm)^{-1}P_n\to (T-\lm)^{-1}$ for every $\lm\in\rho(T)$.

\subsection{Operator matrices with singular Sturm-Liouville operator entries}
We study $2\times 2$ block operator matrices with singular Sturm-Liouville-type operator entries
and prove that the regularization via interval truncation is a spectrally exact approximation.

Let $(a,b)\subset\R$ be a finite interval.
We consider  {Sturm-Liouville differential expressions}\index{operator!Sturm-Liouville} $\tau$ of the form
$$(\tau f)(x):=-(p(x)f'(x))'+q(x)f(x), \quad x\in (a,b),$$
where $p,\,q:(a,b)\to\R$ are measurable functions with $1/p, q \in L^1_{\rm loc}(a,b)$. 
Suppose that there exist $p_{\rm min}>0$ and $q_{\rm min}\in\R$ such that $$p\geq p_{\rm min}, \quad q\geq q_{\rm min}\quad\text{almost everywhere.}$$
We assume that $\tau$ is regular at $b$ and in limit point case at $a$.
Let $(a_n)_{n\in\N}\subset(a,b)$ with 
$a_n\searrow a$, $n\to\infty$. For $n\in\N$ let $P_n:L^2(a,b)\to L^2(a_n,b)$ be the orthogonal projection given by multiplication with the characteristic function of $[a_n,b]$, i.e.\ $P_nf:=\chi_{[a_n,b]}f$, $n\in\N$.
For $\beta\in [0,\pi)$ let $T_{\tau}(\beta)$, $T_{\tau,n}(\beta)$, $n\in\N$, be the selfadjoint realizations of $\tau$ in the Hilbert spaces $L^2(a,b)$, $L^2(a_n,b)$, $n\in\N$, respectively, with domains
\begin{align*}
\dom(T_{\tau}(\beta))&:=\bigg\{f\in L^2(a,b):\,\begin{array}{l}f,pf'\in{\rm AC_{loc}}(a,b),\, \tau f\in L^2(a,b),\\f(b)\cos\beta-(pf')(b)\sin\beta=0\end{array}\bigg\},\\[2mm]
\dom(T_{\tau,n}(\beta))&:=\bigg\{f\in L^2(a_n,b):\,\begin{array}{l}f,pf'\in{\rm AC_{loc}}(a_n,b),\,\tau f\in L^2(a_n,b),\\f(b)\cos\beta-(pf')(b)\sin\beta=0,\,f(a_n)=0\end{array}\bigg\}.
\end{align*}

\begin{theorem}\label{thmSLmatrix}
For $i=1,2$ let $\tau_i$ be a differential expression of the above form, let $\beta_i\in [0,\pi)$ and $\gamma_i\in\C\backslash\{0\}$. Let $s,t, u,v\in L^{\infty}(a,b)$ satisfy $\|s\|_{\infty}\|u\|_{\infty}<|\gamma_1| |\gamma_2|$, and set 
$$s_n:=s|_{[a_n,b]}, \quad t_n:=t|_{[a_n,b]}, \quad u_n:=u|_{[a_n,b]}, \quad v_n:=v|_{[a_n,b]}, \quad n\in\N.$$
Define the orthogonal projections $\cP^{(n)}:={\rm diag}(P_n,P_n)$, $n\in\N$, and the $2\times 2$ block operator matrices $\cA$, $\cA^{(n)}$, $n\in\N$, by
\begin{align*}
\cA&:=\bmat \gamma_1 T_{\tau_1}(\beta_1) & s T_{\tau_2}(\beta_2)+t \\ u T_{\tau_1}(\beta_1)+v & \gamma_2 T_{\tau_2}(\beta_2)\emat, \quad \dom(\cA):=\dom(T_{\tau_1}(\beta_1))\oplus\dom(T_{\tau_2}(\beta_2)),\\
\cA^{(n)}&:=\bmat \gamma_1 T_{\tau_1,n}(\beta_1) & s_n T_{\tau_2,n}(\beta_2)+t_n \\ u_n T_{\tau_1,n}(\beta_1)+v_n & \gamma_2 T_{\tau_2,n}(\beta_2)\emat, \\ \dom(\cA^{(n)})&:=\dom(T_{\tau_1,n}(\beta_1))\oplus\dom(T_{\tau_2,n}(\beta_2)).
\end{align*}
Then there exists $\lm_0\in\C\backslash (\gamma_1 \R\cup\gamma_2\R)$ such that $\cA$, $\cA^{(n)}$, $n\in\N$, satisfy the claims {\rm (a)--(c)}.
If, in addition, $\dom(T_{\tau_1}(\beta_1))=\dom(T_{\tau_2}(\beta_2))$ and  $\dom(T_{\tau_1,n}(\beta_1))=\dom(T_{\tau_2,n}(\beta_2))$, $n\in\N$, then~{\rm~(d)} is satisfied as well.
\end{theorem}

\begin{rem}
The magnetohydrodynamic dynamo operator matrix for the so-called \emph{$\alpha^2$-model} studied in~\cite{dynamopaper} is of similar form as the above matrix. In the more general \emph{$\alpha^2\omega$-model} the operator matrix is an infinite tridiagonal matrix with such $2\times 2$ matrices as building blocks on the diagonal. In order to prove spectral exactness of the interval truncation process for the infinite matrix, we use the results of Subsections~\ref{subsectioninfititegsr} and~\ref{subsectiondisccompmatrices} (see~\cite{dynamopaper}).
\end{rem}

For the proof we need the following lemma.

\begin{lemma}\label{lemmaSL}
Let $\beta\in [0,\pi)$. Then, for every $\lm_0\in\C\backslash\R$, the operators $T_{\tau}(\beta)$, $T_{\tau,n}(\beta)$, $n\in\N$, satisfy the claims  {\rm (a)--(d)}.
\end{lemma}

\begin{proof}
Since the differential expression $\tau$ is in limit point case at the singular endpoint $x=a$, 
the fact that 
$$\Phi_{\tau}(\beta):=\left\{f\in\dom(T_{\tau}(\beta)):\,f=0 \text{ near } x=a\right\}\subset L^2(a,b)$$ 
 is a core of $T_{\tau}(\beta)$ is a well-known result from Sturm-Liouville theory (see e.g.\ the proof of \cite[Satz~14.12]{weid2}). 
For $f\in\Phi_{\tau}(\beta)$ let $n_0(f)\in\N$ be such that $f(x)=0$ for $x\in[a,a_{n_0(f)}]$. This implies $P_nf\in\dom\big(T_{\tau,n}(\beta)\big)$ for $n\geq n_0(f)$.
The strong convergence $P_n\s I, \,n\to\infty$, implies
\beq \label{convforeachxinPhi}
\begin{aligned}
&P_nf\in\dom(T_{\tau,n}(\beta)), \quad n\geq n_0(f),  \\[1mm]
&\big\|T_{\tau,n}(\beta)P_nf-T_{\tau}(\beta)f\big\|
=\big\|P_nT_{\tau}(\beta)f-T_{\tau}(\beta)f\big\|\tolong 0, \quad n\to\infty.
\end{aligned}
\eeq
The selfadjointness of $T_{\tau}(\beta)$, $T_{\tau,n}(\beta)$, $n\in\N$, implies 
\beq\label{eq.nonrealinDeltab}
\C\backslash\R\subset\Delta_b\big((T_{\tau,n}(\beta))_{n\in\N}\big)\cap\rho(T_{\tau}(\beta)).
\eeq
Thus $T_{\tau,n}(\beta)\gsr T_{\tau}(\beta)$ by Theorem~\ref{propcoreforgsr}. Therefore, (c) and (d) are satisfied for every $\lm_0\in\C\backslash\R$.

Now we prove that (b) is satisfied for every $\lm_0\in\C\backslash\R$; the proof of (a) is analogous.
To this end, we show that the embeddings $$B_n: \big(\dom(T_{\tau,n}(\beta)),\|\cdot\|_{T_{\tau,n}(\beta)}\big)\to W^{1,2}(a,b_n), \quad n\in\N,$$
 are uniformly bounded. Then Theorem~\ref{thmsobolev} implies that  the sequence $(J_n)_{n\in\N}$ of embeddings $J_n:\big(\dom(T_{\tau,n}(\beta)),\|\cdot\|_{T_{\tau,n}(\beta)}\big)\to L^2(a,b_n)$, $n\in\N,$ is discretely compact, and hence (b) follows from Proposition~\ref{equivcomp}~ii) and~\eqref{eq.nonrealinDeltab}.
We fix an $n\in\N$ and denote by  $\|\cdot\|_n$, $\langle\cdot,\cdot\rangle_n$ the norm and scalar product of $L^2(a_n,b)$.
Let $f_n\in\dom(T_{\tau,n}(\beta))$ satisfy $\|f_n\|_{T_{\tau,n}(\beta)}=\|f_n\|_{n}+\|T_{\tau,n}(\beta)f_n\|_{n}\leq 1.$
We estimate
\begin{align*}
1&\geq
 \|T_{\tau,n}(\beta)f_n\|_{n}\,\|f_n\|_{n}
\geq \langle T_{\tau,n}(\beta)f_n,f_n\rangle_{n}
 =\int_{a_n}^{b} \big(-(pf_n')'\,\lbar{f_n}+q|f_n|^2\big)(x)\,\rd x\\
 &=\big(-pf_n'\lbar{f_n}\big)(x)\Big|_{x=a_n}^{x=b}+\int_{a_n}^{b} \big(p|f_n'|^2+q|f_n|^2\big)(x)\,\rd x.
 \end{align*}
 If $f_n(b)=0$, then $\big(-pf_n'\lbar{f_n}\big)(x)\big|_{x=a_n}^{x=b}=0$, 
 and if  $(pf_n')(b)=\tan\beta\, f_n(b)$ for some $\beta\in[0,\pi)\backslash\{\pi/2\}$, then
 $\big(-pf_n'\lbar{f_n}\big)(x)\big|_{x=a_n}^{x=b}=-\tan\beta\,|f_n(b)|^2$. If $\beta\in (\pi/2,\pi)$, then the latter is non-negative.
 If $\beta\in  [0,\pi/2)$, then
 \begin{align*}
 \big(-pf_n'\lbar{f_n}\big)(x)\Big|_{x=a_n}^{x=b}
 &=\big(-\tan\beta\,|f_n|^2\big)(x)\Big|_{x=a_n}^{x=b}
 =-\tan\beta\,\int_{a_n}^{b} \frac{\rd}{\rd x}\left(|f_n(x)|^2\right)\,\rd x\\
 &\geq -\tan\beta\,2\|f_n'\|_{n}\,\|f_n\|_{n}
 \geq -\tan\beta\,\left(\eps \|f_n'\|_{n}^2+\frac{1}{\eps}\|f_n\|_{n}^2\right),
 \end{align*} 
 where $\eps>0$ is arbitrary.
 We also use the estimate 
 $$\int_{a_n}^b  \big(p|f_n'|^2\big)(x)\,\rd x\geq p_{\rm min}\|f_n'\|_{n}^2.$$
 Then, if we set $\eps:=p_{\rm min}/(2\tan\beta)$, we obtain altogether 
 \begin{align*}
 1&\geq \frac{1}{2} p_{\rm min}\|f_n'\|_{n}^2+c_{\beta}\|f_n\|_{n}^2,\quad
c_{\beta}:=
 \begin{cases}
 q_{\rm min}, & \beta\in [\pi/2,\pi),\\
 q_{\rm min}-\frac{2 \tan(\beta)^2}{p_{\rm min}}, & \beta\in [0,\pi/2).
 \end{cases}
 \end{align*}
From this it is easy to see that the embeddings $B_n$, $n\in\N$, are uniformly bounded, and thus the claim follows.
\end{proof}

\begin{proof}[Proof of Theorem~{\rm\ref{thmSLmatrix}}]
By the assumption $\|s\|_{\infty}\|u\|_{\infty}<|\gamma_1| |\gamma_2|$,
every sufficiently small $\eps>0$ satisfies
\beq\label{eq.epschoice}
\left(\frac{\|u\|_{\infty}}{|\gamma_1|}(1+\eps)+\|v\|_{\infty}\eps\right)\left(\frac{\|s\|_{\infty}}{|\gamma_2|}(1+\eps)+\|t\|_{\infty}\eps\right)<1.
\eeq
For such an $\eps$, choose $\lm_0\in\C\backslash (\gamma_1 \R\cup\gamma_2\R)$ so that
$$\sup_{\xi\in\gamma_i\R}\frac{|\xi|}{|\xi-\lm_0|}\leq 1+\eps, \quad {\rm dist}(\lm_0,\gamma_i\R)\geq \frac{1}{\eps},\quad i=1,2.$$
Then, using the selfadjointness of $T_{\tau,n}(\beta_i)$, $i=1,2$, and~\cite[Equation~V.(3.17)]{kato}, 
\begin{align*}
\|(u_nT_{\tau,n}(\beta_1)+v_n)(\gamma_1 T_{\tau,n}(\beta_1)-\lm_0)^{-1}\|&\leq \frac{\|u\|_{\infty}}{|\gamma_1|}\sup_{\xi\in\gamma_1\R}\frac{|\xi|}{|\xi-\lm_0|}+\frac{\|v\|_{\infty}}{{\rm dist}(\lm_0,\gamma_1\R)}\\
&\leq \frac{\|u\|_{\infty}}{|\gamma_1|}(1+\eps)+\|v\|_{\infty}\eps=:\gamma_{\lm_0}^{(1)},\\
\|(s_nT_{\tau,n}(\beta_1)+t_n)(\gamma_2 T_{\tau,n}(\beta_2)-\lm_0)^{-1}\|&\leq \frac{\|s\|_{\infty}}{|\gamma_2|}\sup_{\xi\in\gamma_2\R}\frac{|\xi|}{|\xi-\lm_0|}+\frac{\|t\|_{\infty}}{{\rm dist}(\lm_0,\gamma_2\R)}\\
&\leq \frac{\|s\|_{\infty}}{|\gamma_2|}(1+\eps)+\|t\|_{\infty}\eps=:\gamma_{\lm_0}^{(2)}.
\end{align*}
Note that $\gamma_{\lm_0}^{(1)} \gamma_{\lm_0}^{(2)}<1$ by~\eqref{eq.epschoice}.
Then Lemma~\ref{lemmaSL} and Theorem~\ref{propmatrixwithdiscrcompres}  imply that $\big((\cA^{(n)}-\lm_0)^{-1}\big)_{n\in\N}$ is discretely compact, i.e.\ (b) is satisfied. 
Claim (a) is shown analogously.

The generalized strong resolvent convergence in (c) follows from Lemma~\ref{lemmaSL}, \eqref{convforeachxinPhi} in its proof, and from Theorem~\ref{propgsrformatrix} using  that $\lm_0\in\Delta_b\big((\cA^{(n)})_{n\in\N}\big)\cap\rho(\cA)\neq\emptyset$. 

To prove (d), we first note that the operator matrices $\cA$, $\cA^{(n)}$, $n\in\N$, are relatively bounded perturbations of diagonal operators with relative bound $<1$ (by~\eqref{eq.epschoice} and Remark~\ref{remgammaorder}), and the same holds for the adjoint matrices if we assume that, in addition, $\dom(T_{\tau_1}(\beta_1))=\dom(T_{\tau_2}(\beta_2))$ and  $\dom(T_{\tau_1,n}(\beta_1))=\dom(T_{\tau_2,n}(\beta_2))$ for $n\in\N$. Then~\cite[Corollary~1]{Hess-Kato} yields 
\begin{alignat*}{3}
\cA^*&=\bmat \overline{\gamma_1} T_{\tau_1}(\beta_1) & \overline{u} T_{\tau_1}(\beta_1)+\overline{v}\\  \overline{s} T_{\tau_2}(\beta_2)+\overline{t} & \overline{\gamma_2} T_{\tau_2}(\beta_2)\emat, \quad& \dom(\cA^*)&=\dom(\cA),\\
\big(\cA^{(n)}\big)^*&=\bmat \overline{\gamma_1} T_{\tau_1,n}(\beta_1) & \overline{u_n} T_{\tau_1,n}(\beta_1)+\overline{v_n}\\ \overline{s_n} T_{\tau_2,n}(\beta_2)+\overline{t_n} & \overline{\gamma_2} T_{\tau_2,n}(\beta_2)\emat, \quad&
\dom\big(\big(\cA^{(n)}\big)^*\big)&=\dom\big(\cA^{(n)}\big).
\end{alignat*}
Now the strong convergence $\big(\big(\cA^{(n)}\big)^*-\overline{\lm_0}\big)^{-1}\cP^{(n)}\s (\cA^*-\overline{\lm_0})^{-1}$ is shown analogously as (c).
\end{proof}

\subsection{Domain truncation of magnetic Schr\"odinger operators on $\R^d$}
Let $\Omega_n\subset \R^d$, $n\in\N$, be nested open sets that exhaust $\R^d$ eventually.
Denote by $\|\cdot\|$, $\langle\cdot,\cdot\rangle$, $\|\cdot\|_n$, $\langle\cdot,\cdot\rangle_n$, $n\in\N$, the norm and scalar product of $L^2(\R^d)$ and $L^2(\Omega_n)$, $n\in\N$, respectively.
Let $P_n$, $n\in\N$, be the orthogonal projections of $L^2(\R^d)$ onto the respective subspaces, 
given by multiplication with the characteristic function~$\chi_{\Omega_n}$.
Then $P_n\s I$ as $n\to\infty$.

Consider the differential expression $\tau:=-\Delta+p\cdot\nabla+v$ with a vector potential $p:\R^d\to\C^d$ and a scalar potential $v:\R^d\to\C$.
An important application is given by the \emph{magnetic Schr\"odinger operator} $\tau=(-\I\nabla+ A)^2+V=-\Delta+p\cdot\nabla+v$ with $p=-\I A$ and $v=-\I\nabla\cdot A+A^2+V$.

We assume that $p$ and $v=q+r$ satisfy 
\begin{enumerate}
\item[\rm (i)]  $p\in L^{\infty}(\R^d)$;
\item[\rm (ii)] $q\in{\rm AC_{loc}}(\R^d)$, $\re\,q\geq 0$, $|q(x)|\to\infty$ as $|x|\to\infty$, and there are $a_{\nabla},b_{\nabla}\geq 0$ with 
\beq \label{eq.nabla}|\nabla q(x)|^2\leq a_{\nabla} + b_{\nabla}\, |q(x)|^2\quad \text{for almost all } x\in\R^d;\eeq
\item[\rm (iii)] $r\in L^2_{\rm loc}(\R^d)$, and
there exist $a_r,b_r\geq 0$ with $b_r<1$ such that
\beq |r(x)|^2\leq a_r + b_r |q(x)|^2 \quad\text{for almost all } x\in\R^d.\label{eq.relbddpot}\eeq
\end{enumerate}

\begin{theorem}\label{thm.schroed}
Let $A$ and $A_n$, $n\in\N$, be the Dirichlet realizations of $\tau $ in the respective spaces,
\begin{alignat*}{3}
Af&:=\tau f, \quad& \dom(A)&:=
\big\{f\in W^{2,2}(\R^d):\,qf\in L^2(\R^d)\big\},\\
A_nf&:=\tau f, \quad& \dom(A_n)&:=
W^{2,2}(\Omega_n)\cap W_0^{1,2}(\Omega_n), \quad n\in\N.
\end{alignat*}
Then, for every real $\lm_0<0$ with sufficiently large $|\lm_0|$, the operators $A$, $A_n$, $n\in\N$, satisfy the claims  {\rm (a)--(c)}.
If, in addition, $p\in W^{1,\infty}(\R^d)$, then~{\rm~(d)} is satisfied as well.
\end{theorem}

For the proof we use the following result.

\begin{lemma}\label{lemma.schroed}
Let $T$ and $T_n$, $n\in\N$, be the Dirichlet realizations of  the Schr\"odinger differential expression $\tau_0:=-\Delta+q$ in the respective spaces,
\begin{alignat*}{3}
 Tf&:=\tau_0 f, \quad& \dom(T)&:=
\big\{f\in W^{2,2}(\R^d):\,qf\in L^2(\R^d)\big\},\\
 T_nf&:=\tau_0 f, \quad& \dom(T_n)&:=
W^{2,2}(\Omega_n)\cap W_0^{1,2}(\Omega_n), \quad n\in\N.
\end{alignat*}
Then,  for every real $\lm_0<0$, the operators $T$, $T_n$, $n\in\N$, satisfy the claims  {\rm (a)--(d)}.
\end{lemma}

\begin{proof}
To prove claim (b), we fix an $n\in\N$. First note that $q_n:=q|_{\Omega_n}$ is bounded and thus $T_n$ is a bounded perturbation of the Dirichlet Laplacian on $\Omega_n$.
By~\cite[Theorem~VI.1.4]{edmundsevans}, the operator $T_n$ is $m$-accretive with compact resolvent. 
Therefore, every $\lm_0<0$ satisfies $\lm_0\in\rho(T_n)$ and, using \cite[Problem~V.3.31]{kato},
\beq \|(T_n-\lm_0)^{-1}\|\leq \frac{1}{|\lm_0|}, \quad \|T_n(T_n-\lm_0)^{-1}\|\leq 1.\label{eq.accretive}\eeq
So we have, in particular, 
\beq\label{eq.neglminDeltab}
\forall\,\lm_0<0:\quad \lm_0 \in\Delta_b\big((T_n)_{n\in\N}\big).
\eeq

For any $f\in \dom(T_n)$ we obtain, using integration by parts and $\re\,q_n\geq 0$ by assumption~(ii), 
$$\|T_nf\|_n^2=\!\|\Delta f\|_n^2+\|q_nf\|_n^2+2\re \langle -\Delta f,q_nf\rangle_n\!\geq\! \|\Delta f\|_n^2+\|q_nf\|_n^2+2\re \langle \nabla f,\!(\nabla q_n)f\rangle_n.$$ 
Now, again using integration by parts  and with~\eqref{eq.nabla}, for any $\eps,\delta>0$,
\begin{align*}
&\big| 2 \re \langle \nabla f,(\nabla q_n)f\rangle_n\big|\\
&\leq \frac{1}{\eps}\|\nabla f\|_n^2+\eps \|(\nabla q_n)f\|_n^2
\leq \frac{1}{\eps}\langle -\Delta f,f\rangle_n + \eps a_{\nabla}\|f\|_n^2+\eps b_{\nabla} \|q_nf\|_n^2\\
&\leq \bigg(\eps a_{\nabla}+\frac{1}{4 \eps\delta}\bigg)\|f\|_n^2+\frac{\delta}{\eps}\|\Delta f\|_n^2+\eps b_{\nabla}\|q_nf\|_n^2.
\end{align*}
Let $\alpha\in (0,1)$ be arbitrary. We choose $\eps>0$ and then $\delta>0$ both so small that $\max\{\delta/\eps,\eps b_{\nabla}\}\leq \alpha$.
With these $\eps, \delta$, we set $C_{\alpha}:=\eps a_{\nabla} +1/(4\eps\delta)$  and arrive at 
\beq \|T_nf\|_n^2+C_{\alpha} \|f\|_n^2\geq (1-\alpha)\big(\|\Delta f\|_n^2+\|q_nf\|_n^2\big).\label{eq.normequiv}\eeq

Now let $I\subset\N$ be an infinite subset and let $f_n\in\dom(T_n)$, $n\in I$, be such that 
the sequence of graph norms $(\|f_n\|_{T_n})_{n\in I}$ is  bounded. Then~\eqref{eq.normequiv} implies that $(\|\Delta f_n\|_n)_{n\in I}$ and $(\|q_nf_n\|_n)_{n\in I}$ are bounded, and hence so is $(\|\nabla f_n\|_n)_{n\in I}$ since $2\|\nabla f_n\|_n^2\leq \|f_n\|_n^2+\|\Delta f_n\|_n^2$, $n\in I.$
By extending every $f_n$ by zero outside its domain~$\Omega_n$, we obtain $(f_n)_{n\in I}\subset W^{1,2}(\R)$, and 
the sequences $(\|\nabla f_n\|)_{n\in I}$ and 
$(\|qf\|)_{n\in I}$ are both bounded. 
Since $|q(x)|\to\infty$ as $|x|\to\infty$ by assumption~(ii), \emph{Rellich's criterion} (see \cite[Theorem~XIII.65]{reedsimon}) implies that $(f_n)_{n\in I}$ has a subsequence that is convergent in $L^2(\R)$.
 So we have shown that the sequence  $(J_n)_{n\in\N}$ of embeddings $J_n:\big(\dom(T_{n}),\|\cdot\|_{T_{n}}\big)\to L^2(\Omega_n)$, $n\in\N,$ is discretely compact, and hence, by Proposition~\ref{equivcomp}~ii) and~\eqref{eq.neglminDeltab}, claim (b) follows for every $\lm_0<0$.

Let $T_{\min}, T_{\max}$ be the realizations of $\tau_0=-\Delta+q$ with 
$$\dom(T_{\min})=\Phi:=C_0^{\infty}(\R^d), \quad \dom(T_{\max}):=\{f\in L^2(\R^d):\,\tau f\in L^2(\R^d)\}.$$
By~\cite[Theorem~VII.2.6, Corollary~VII.2.7]{edmundsevans}, $T_{\max}$ is $m$-accretive and the closure of $T_{\min}$.
Be proceeding as before, we arrive at the inequality~\eqref{eq.normequiv} with $T_n$ replaced by $T_{\min}$;
since $T_{\max}=\overline{T_{\min}}$, the inequality also holds for $T_{\max}$.
Therefore, it is easy to see that $T_{\max}=T$. So $T$ is $m$-accretive and $\Phi$ is a core of $T$.
Claim~(a) is shown analogously as (b).

For every $f\in\Phi$ there exists $n_0(f)\in\N$ such that $f\in\dom(T_n)$ and $T_nf=P_nTf$ for $n\geq n_0(f)$. 
Hence Theorem~\ref{propcoreforgsr} and~\eqref{eq.neglminDeltab} imply that (c) is satisfied for every $\lm_0<0$.
Since, by \cite[Theorem~VII.2.5]{edmundsevans}, the adjoint operators $T^*$, $T_n^*$, $n\in\N$, are the Dirichlet realizations of $\tau^*=-\Delta+\overline{q}$ and $\overline{q}$ satisfies the assumption~(ii), we obtain analogously that (d) is satisfied for every $\lm_0<0$.
\end{proof}

\begin{proof}[Proof of Theorem~{\rm\ref{thm.schroed}}]
Let $\lm_0<0$. The assumption~(i) yields, for any~$\beta>0$,
\beq\label{eq.pnabla}
\| p\cdot \nabla f\|^2\leq \|p\|_{\infty}^2\|\nabla f\|^2\leq \frac{\|p\|_{\infty}^4}{4\beta}\|f\|^2+\beta \|\Delta f\|^2.
\eeq
Define $p_n:=p|_{\Omega_n}$  and $r_n:=r|_{\Omega_n}$ for $n\in\N$, and $S:=A-T$, $S_n:=A_n-T_n$, $n\in\N$.
 Then the estimates~\eqref{eq.pnabla} and~\eqref{eq.relbddpot} imply, for any $\nu>0$,
\begin{align*}
\|S_n(T_n-\lm_0)^{-1}\|^2
&\leq \left(1+\frac{1}{4\nu}\right)\|(p\cdot \nabla)(T_n-\lm_0)^{-1}\|^2+(1+\nu)\|r_n(T_n-\lm_0)^{-1}\|^2\\
&\leq \left(\frac{\|p\|_{\infty}^4}{4\beta}\left(1+\frac{1}{4\nu}\right)+a_r(1+\nu)\right)\|(T_n-\lm_0)^{-1}\|^2\\
&\quad+\beta\left(1+\frac{1}{4\nu}\right)\|\Delta(T_n-\lm_0)^{-1}\|^2+ b_r(1+\nu)\|q_n(T_n-\lm_0)^{-1}\|^2.
\end{align*}
We choose $\nu$ and then $\beta$ so small that $b:=\max\{\beta(1+1/(4\nu)),b_r(1+\nu)\}<1$.
Set $$a:=\frac{\|p\|_{\infty}^4}{4\beta}\left(1+\frac{1}{4\nu}\right)+a_r(1+\nu).$$
Then, by~\eqref{eq.normequiv} and~\eqref{eq.accretive},
\begin{align*}
\|S_n(T_n-\lm_0)^{-1}\|^2
&\leq  \bigg(a+\frac{ b C_{\alpha}}{1-\alpha}\bigg)\|(T_n-\lm_0)^{-1}\|^2+\frac{b}{1-\alpha}\|T_n(T_n-\lm_0)^{-1}\|^2\\
& \leq  \bigg(a+\frac{ b C_{\alpha}}{1-\alpha}\bigg)\frac{1}{|\lm_0|^2}+\frac{b}{1-\alpha}=:(\gamma_{\lm_0,\alpha})^2.
\end{align*}
Now we choose $\alpha\in (0,1)$ so small that $b/(1-\alpha)<1$ and then $|\lm_0|$ so large that $\gamma_{\lm_0,\alpha}<1$. Then 
Theorem~\ref{perturbdiscrcompres} yields that  $\lm_0\in\rho(A_n)$, $n\in\N$, and $\big((A_n-\lm_0)^{-1}\big)_{n\in\N}$ is discretely compact, i.e.\ claim~(b)  holds for this $\lm_0$. Claim~(a) is shown analogously.

Analogously as $\|S_n(T_n-\lm_0)^{-1}\|\leq\gamma_{\lm_0,\alpha}<1$, $n\in\N$, if $|\lm_0|$ is sufficiently large,  one may show that $\|S(T-\lm_0)^{-1}\|<1$ if $|\lm_0|$ is sufficiently large. 
Now claim~(c) follows from Theorem~\ref{perturbgsr} provided that $S_n(T_n-\lm_0)^{-1}P_n\s S(T-\lm_0)^{-1}$. 
To show the latter, we take $f\in\Phi=\C_0^{\infty}(\R^d)$; the latter is a core of $T$ (see the proof of Lemma~\ref{lemma.schroed}).
Define $g:=(T-\lm_0)f$. Then there exists $n_0(f)\in\N$ such that ${\rm supp}\,g\subset{\rm supp}\,f\subset \Omega_n$, $n\geq n_0(f)$.
Hence, for $n\geq n_0(f)$, $$S_n(T_n-\lm_0)^{-1}g=S_n(T_n-\lm_0)^{-1}(T-\lm_0)f=S_nf=Sf=S(T-\lm_0)^{-1}g.$$
Now $S_n(T_n-\lm_0)^{-1}P_n\s S(T-\lm_0)^{-1}$ follows since $\{(T-\lm_0)f:\,f\in\Phi\}\subset L^2(\R^d)$ is a dense subset and $\|S_n(T_n-\lm_0)^{-1}\|$, $n\in\N$, are uniformly bounded.

Now assume that, in addition, $p\in W^{1,\infty}(\R^d)$. Then $A^*$ and $A_n^*$, $n\in\N$, are Dirichlet realizations of the adjoint differential expression
$$\tau^*=-\Delta-\overline{p}\cdot\nabla+(-\nabla\cdot \overline{p}+\overline{v}).$$
Since the vector potential $\widetilde p:=\overline{p}$ and the scalar potential $\widetilde v:=\widetilde q+\widetilde r$ with $\widetilde q:=\overline{q}$ and $\widetilde r:=\overline{r}-\nabla\cdot\overline{p}$ satisfy assumptions~(i)--(iii), the above arguments imply that $(A_n^*-\lm_0)^{-1}P_n\s (A^*-\lm_0)^{-1}$ for every $\lm_0<0$ with $|\lm_0|$ large enough; thus (d) is satisfied.
\end{proof}

\subsection{Galerkin approximation of block-diagonally dominant matrices}

Let $\{e_k:\,k\in\N\}$ be the standard orthonormal basis of $H:=l^2(\N)$. Define
the $k$-dimensional subspace $H_k:={\rm span}\{e_i:\,i=1,\dots,k\}\subset H$.
Denote by $P_k:H\to H_k$, $k\in\N$, the corresponding orthogonal projections. Obviously, $P_k\s I$ as $k\to\infty$.

We study the Galerkin approximation of a closed operator $A\in C(H)$. To this end, we identify $A$ with its matrix representation with respect to $\{e_k:\,k\in\N\}$,
$$A=(A_{ij})_{i,j=1}^{\infty}, \quad A_{ij}=\langle A e_j,e_i\rangle, \quad i,j\in\N.$$
With $k_0:=0$ and a strictly increasing sequence $(k_n)_{n\in\N}\subset\N$ define the diagonal blocks
$$B_n:=(A_{ij})_{i,j=k_{n-1}+1}^{k_n}, \quad n\in\N,$$
and split $A$ as $A=T+S$ with $T:={\rm diag}(B_n:\,n\in\N)$.
Define the Galerkin approximations $A_n:=P_{k_n}A|_{H_{k_n}}$, $n\in\N$.

\begin{theorem}\label{thmgalerkin}
Assume that there exists $\lm_0\in\underset{n\in\N}{\bigcap}\rho(B_n)$ with
$\|(B_n-\lm_0)^{-1}\|\to 0$ as $n\to\infty$. Then $\lm_0\in\rho(T)$.
If $\dom(T)\subset\dom(S)$ and $\|S(T-\lm_0)^{-1}\|<1$,
then $A$, $A_n$, $n\in\N$, satisfy the claims  {\rm (a)--(c)}.
If, in addition,   $\dom(T^*)\subset\dom(S^*)$ and $\|S^*(T^*-\overline{\lm_0})^{-1}\|<1$, then~{\rm~(d)} is satisfied as well.
\end{theorem}

The proof relies on the following lemma.

\begin{lemma}\label{lemmablockdiag}
If $\|(B_n-\lm_0)^{-1}\|\to 0$ as $n\to\infty$, then the block-diagonal operators $T$ and $T_n:={\rm diag}(B_k:\,k=1,\dots,n)$, $n\in\N$, satisfy the claims  {\rm (a)--(d)}.
\end{lemma}

\begin{proof}
The assumption $\|(B_n-\lm_0)^{-1}\|\to 0$ implies that $(T-\lm_0)^{-1}$ is the norm limit of the finite-rank (and thus compact) operators $(T_n-\lm_0)^{-1}P_{k_n}$ and hence compact by \cite[Theorem III.4.7]{kato}.
In addition, $(T_n-\lm_0)^{-1}=P_{k_n}(T-\lm_0)^{-1}|_{H_{k_n}}$ and therefore $(T_n-\lm_0)^{-1}$, $n\in\N$, are compact and form a discretely compact sequence.
Thus (a) and (b) are satisfied. 

Again using $\|(B_n-\lm_0)^{-1}\|\to 0$, we see that $(T_n-\lm_0)^{-1}P_{k_n}\to (T-\lm_0)^{-1}$ which implies, in particular, that (c) and (d) hold.
\end{proof}

\begin{proof}[Proof of Theorem~{\rm\ref{thmgalerkin}}]
Define the Galerkin approximations $S_n:=P_{k_n}S|_{H_{k_n}}$, $n\in\N$.
Note that $$S_n(T_n-\lm_0)^{-1}=P_{k_n}S(T-\lm_0)^{-1}|_{H_{k_n}},\quad n\in\N.$$
So we readily conclude that $\|S_n(T_n-\lm_0)^{-1}\|\leq \|S(T-\lm_0)^{-1}\|=:\gamma_{\lm_0}<1$ and $S_n(T_n-\lm_0)^{-1}P_{k_n}\s S(T-\lm_0)^{-1}$.
Using Lemma~\ref{lemmablockdiag} and 
Theorems~\ref{perturbgsr},~\ref{perturbdiscrcompres}, we obtain claims~(d) and (c)  (and (a) analogously).

If the additional assumptions $\dom(T^*)\subset\dom(S^*)$ and $\|S^*(T^*-\overline{\lm_0})^{-1}\|<1$ hold, then claim~(d) is proved analogously; note that $(T+S)^*=T^*+S^*$ by~\cite[Corollary~1]{Hess-Kato}.
\end{proof}

\begin{rem}
If the assumptions of Theorem~\ref{thmgalerkin} are satisfied, then the Galerkin approximation $(A_n)_{n\in\N}$ with $A_n:=P_{k_n}A|_{H_{k_n}}$, $n\in\N$, is spectrally exact.
However, if we consider all $P_kA|_{H_k}$, $k\in\N$, then spurious eigenvalues may occur.
As an example, let $A$ be the selfadjoint Jacobi operator
$$A=\bmat 0 & q_1  & 0 & \hdots \\ q_1 & 0 & q_2 & \ddots \\ 0 & q_2 & 0 &\ddots & \\ \vdots &\ddots &\ddots  & \ddots & \emat, \quad 
q_k:=\begin{cases} k+1, & k\text{ odd},\\ k/2, & k\text{ even}.\end{cases}$$
One may check that the assumptions of Theorem~\ref{thmgalerkin} are satisfied for $\lm_0=0$ and
$$k_n=2n, \quad B_n=\bmat 0 & q_{2n-1} \\ q_{2n-1} & 0 \emat, \quad n\in\N.$$
So the operators $P_{2n}A|_{H_{2n}}$, $n\in\N$, form a spectrally exact approximation of $A$.
However, the point $\lm_0=0\in\rho(A)$ is an eigenvalue of every $P_{2n-1}A|_{H_{2n-1}}$, $n\in\N$, and thus a point of spectral pollution.
\end{rem}

\subsection*{Acknowledgements}
This paper is based on Chapter~1 of the author's Ph.D.\ thesis~\cite{boegli-phd}. She would like to thank her doctoral advisor Christiane Tretter for the guidance.
The work was supported by 
the Swiss National Science Foundation (SNF), grant no.\ 200020\_146477 and Early Postdoc.Mobility project P2BEP2\_159007.

{\footnotesize
\bibliographystyle{acm}
\bibliography{mybib}
}

\end{document}